\documentclass{nmeauth}

\usepackage[pdftex,colorlinks,bookmarksopen,bookmarksnumbered,citecolor=red,urlcolor=red]{hyperref}

\usepackage{graphicx,color}
\usepackage{amsmath,amssymb,amsfonts}
\usepackage{fancybox}
\usepackage{subfigure}
\usepackage{float}
\usepackage{multirow}
\usepackage{bbdd} 
\usepackage{xspace}
\usepackage{mynotations}

\theoremstyle{plain}

\newtheorem{remark}{\itshape \rmfamily Remark}

\newcommand{\Figure}{Figure}
\newcommand{\Figures}{Figures}
\newcommand{\Table}{Table}

\newcommand{\Section}{Section}

\newcommand{\Fig}[1]{\Figure~\ref{#1}}
\newcommand{\Figs}[1]{\Figures~\ref{#1}}
\newcommand{\Tab}[1]{\Table~\ref{#1}}

\newcommand{\Sect}[1]{\Section~\ref{#1}}

\newcommand{\ie}{i.e.\xspace}

\begin{document}
\NME{1}{6}{00}{28}{10}

\runningheads{F.\ Pled, L.\ Chamoin and P.\ Ladev\`eze}
{On techniques for constructing admissible stress fields in model verification}

\title{On the techniques for constructing admissible stress fields in model verification: performances on engineering examples}
\author{F.~Pled\affil{1}, L.~Chamoin\affil{1}, \ and P.~Ladev\`eze\affil{1}\comma\affil{2}\comma\corrauth}

\address{\affilnum{1}\ LMT-Cachan (ENS-Cachan/CNRS/Paris 6 University), 61 Avenue du Pr\'esident Wilson, 94235 CACHAN Cedex, FRANCE\\
\affilnum{2}EADS Foundation Chair, Advanced Computational Structural Mechanics, FRANCE}

\corraddr{P. Ladev\`eze, LMT-Cachan (ENS-Cachan/CNRS/Paris 6 University), 61 Avenue du Pr\'esident Wilson, 94235 CACHAN Cedex}



\received{30 September 2010}
\revised{27 January 2011}
\accepted{4 February 2011}


\begin{abstract}
Robust global/goal-oriented error estimation is used nowadays to control the approximate finite element solutions obtained from simulation. In the context of Computational Mechanics, the construction of admissible stress fields (\ie stress tensors which verify the equilibrium equations) is required to set up strict and guaranteed error bounds (using residual based error estimators) and plays an important role in the quality of the error estimates. This work focuses on the different procedures used in the calculation of admissible stress fields, which is a crucial and technically complicated point. The three main techniques that currently exist, called the element equilibration technique (EET), the star-patch equilibration technique (SPET), and the element equilibration + star-patch technique (EESPT), are investigated and compared with respect to three different criteria, namely the quality of associated error estimators, computational cost and easiness of practical implementation into commercial finite element codes.

The numerical results which are presented focus on industrial problems; they highlight the main advantages and drawbacks of the different methods and show that the behavior of the three estimators, which have the same convergence rate as the exact global error, is consistent. Two- and three-dimensional experiments have been carried out in order to compare the performance and the computational cost of the three different approaches. The analysis of the results reveals that the SPET is more accurate than EET and EESPT methods, but the corresponding computational cost is higher. Overall, the numerical tests prove the interest of the hybrid method EESPT and show that it is a correct compromise between quality of the error estimate, practical implementation and computational cost. Furthermore, the influence of the cost function involved in the EET and the EESPT is studied in order to optimize the estimators.
\end{abstract}


\section{INTRODUCTION}\label{1}
Assessing the global/goal-oriented discretization error, henceforth known as \textit{model verification}, has become a major challenge and a topical issue for both industrial applications and academic research. The widespread availability of computer hardware and numerical tools has contributed to a recent upsurge in the development of virtual prototyping. The modeling of physical problems requires the use of an initial mathematical model, which is considered as the reference to build a discretized model whose calculation is performed by numerical methods suited to computing tools. One of the most widespread approximation methods is the Finite Element Method (FEM) providing a numerical approximated solution of the \textit{a priori} unknown exact solution of the reference model.

Within the finite element framework, the first works dealing with verification date back from the late 1970s \cite{Lad75,Lad83,Bab78a,Zie87} and provide a global estimation of the discretization error which allows the global quality of a FE calculation to be quantified \cite{Ver96,Bab01,Ste03,Lad04}. For over fifteen years, various methods had been introduced to control the numerical quality of specific quantities of interest, such as local stresses, displacement values; they lead to bounds on functional outputs which are relevant information for design purposes in Mechanics \cite{Par97,Ran97,Cir98,Per98,Lad99a,Pru99,Str00,Bec01,Wib06,Lad08,Cha07,Cha08,Pan10}. However, among those global/goal-oriented error estimators built from different methods, few lead to guaranteed and relevant bounds, which is a serious drawback in the domain of robust design. Only methods based on the construction of statically admissible stress fields actually yield guaranteed bounds of the global discretization error and allow the treatment of a wide range of mechanical problems (visco-elasticity, visco-plasticity, transient dynamics, ...).\\

Various techniques have been developed for the construction of such stress fields; the first technique, proposed by Fraeijs de Veubeke, is based on a dual formulation, which corresponds to the best approach as regards the quality of the error estimator \cite{Fra65,Kem03}. However, despite its significant performances, this method is very expensive and intractable in practice as it requires the calculation of another global approximated solution of the reference problem; therefore its implementation is not suited for current FE codes. Three other techniques are well-suited to error estimation as they provide a statically admissible stress field from the FE solution: the element equilibration technique (EET) \cite{Lad83,Ain93bis,Lad96,Lad97,Flo02}; the star-patch equilibration technique (SPET) \cite{Mac00,Car00,Mor03,Pru04,Par06,Par09,Moi09,Cot09,Gal09}; the element equilibration + star-patch technique (EESPT) \cite{Lad08bis,Lad10bis}.

First, the EET, also called \textit{standard technique}, introduces equilibrated tractions along finite element edges, leading to a quasi-explicit calculation of admissible stress fields at the element scale. The construction of tractions, in equilibrium with the external loading, is based on properties of the finite element stress field through an energy relation called the \textit{prolongation condition}. Despite its affordable computational cost, the procedure for defining such tractions is in general difficult to implement into existing finite element codes. On the contrary, the SPET, also called the \textit{flux-free technique}, precludes the main drawbacks of the EET as it circumvents the need of flux-equilibration. The principle of the SPET is to exploit a partition of unity in order to define self-equilibrated local problems (it avoids the calculation of equilibrated tractions) and results in a simple implementation. However, this technique requires the fine solution of local problems on sets of elements, also called \textit{patches} or \textit{stars}, involving a large number of degrees of freedom (particularly for 3D applications) and, therefore, leading to the calculation of admissible stress fields at high computational cost. Eventually, the EESPT is a hybrid technique which takes advantage of the ingredients of both EET and SPET methods. Indeed, this technique combines the prolongation condition used in the EET and the partition of unity involved in the SPET in order to construct equilibrated tractions over element edges. This hybrid procedure comes down to solving simple, local, independent problems defined over patches of elements, then solving problems at the element scale at reasonable computational cost.\\

This paper aims at comparing these three techniques for constructing statically admissible stress fields in terms of quality of the associated error estimate, practical implementation and computational cost. The performance of error estimators is investigated by comparing the estimate with the energy norm of the exact error. The numerical experiments are carried out on industrial cases. Thus, this work is in line with previous studies even though it goes further than that presented in \cite{Lad10bis}, which compared the three error estimation methods on two-dimensional academic examples, and provides more details on implementation. A good match is observed between the estimates provided by the three techniques and the energy norm of the reference error. Furthermore, different cost functions involved in the minimization step of both EET and EESPT methods are considered and analyzed to assess their influence on the corresponding error estimates.

The remainder of the paper is structured as follows: after this introduction, \Sect{2} defines the reference problem and its discretized form in order to introduce the basic concept of admissible solution; \Sect{3} deals with the main principles of the three existing methods for constructing admissible stress fields, while \Sect{4} is devoted to a detailed analysis of their practical implementation; several two- and three-dimensional numerical results are presented in \Sect{5} to investigate the advantages and drawbacks of those techniques on complex industrial structures; finally, \Sect{6} suggests several perspectives on which future research should focus.

\section{BASICS ON ERROR ESTIMATION AND ADMISSIBLE SOLUTIONS}\label{2}

\subsection{Presentation of the reference problem}\label{2.1}

\begin{figure}
\centering\includegraphics[scale=0.22]{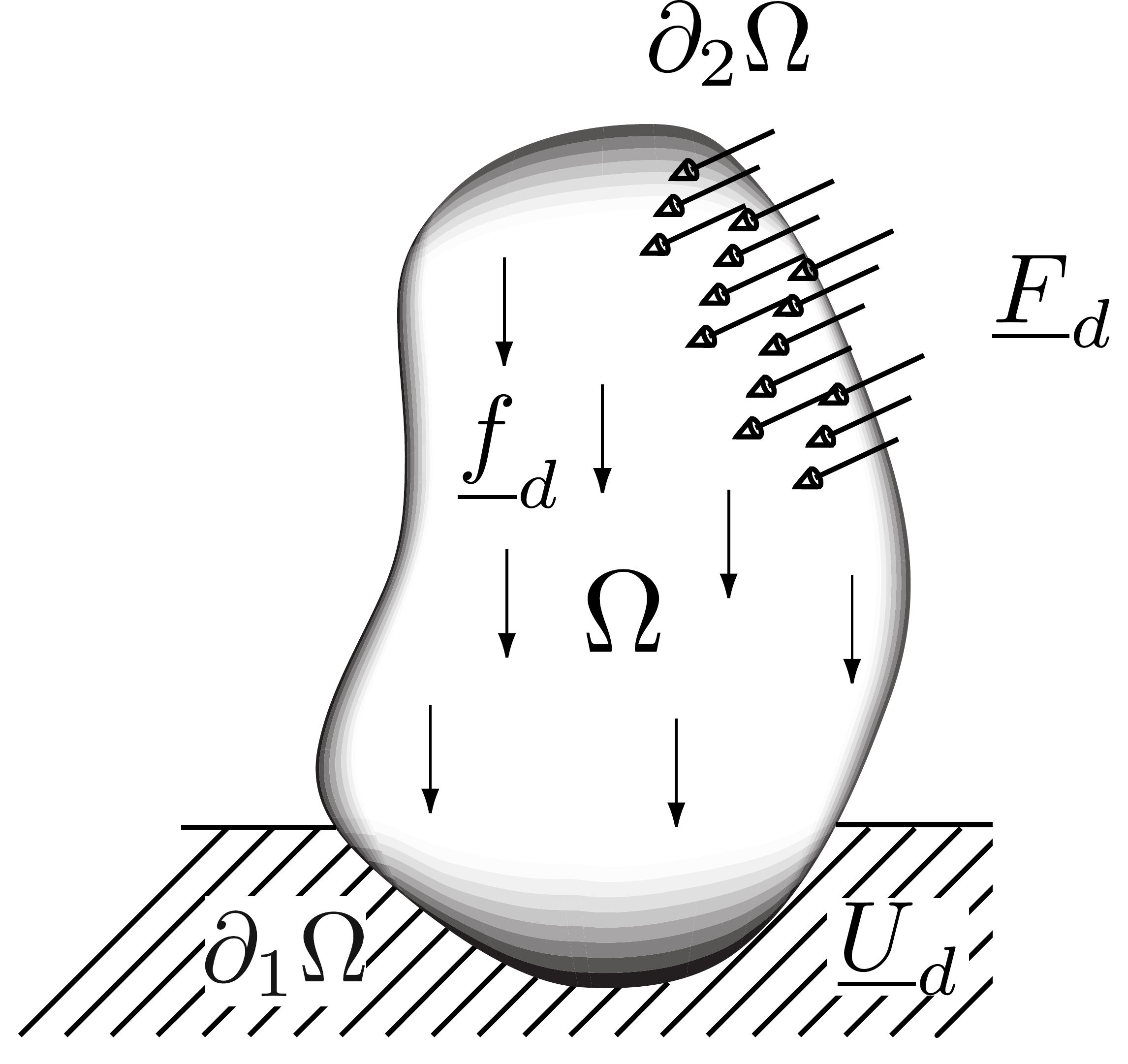}
\caption{The structure and its environment.}\label{fig1:pbref_black_wo_line}
\end{figure}

Let us consider an open bounded domain $\Om$. Its boundary $\dOm$ is divided into two disjoints, complementary parts $\partial_1 \Om\, (\neq \emptyset)$ and $\partial_2 \Om$ such that $\overline{\partial_1 \Om \cup \partial_2 \Om} = \dOm$, $\partial_1 \Om \cap \partial_2 \Om = \emptyset$. The problem is assumed to be linear under the assumptions of elastic material, small perturbations and the evolution in time is considered to be quasi-static and isothermal. The mechanical structure is subjected to a loading represented by (see \Fig{fig1:pbref_black_wo_line}):

\begin{itemize}
\item
a prescribed displacement $\Uu_d$ on part $\partial_1 \Om$;\\
\item
a traction force density $\und{F}_d$ on part $\partial_2 \Om$;\\
\item
a body force field $\und{f}_d$ in $\Om$.
\end{itemize}

The reference problem consists of searching a displacement/stress pair $(\uu (\Mu),\cont (M)), M \in \Om$, which verifies:
\begin{itemize}
\item 
the kinematic conditions:
\begin{multline}\label{eq1:CAref}
\shoveright{\uu \in \Ucb; \quad \uu_{\restrictto{\partial_1 \Om}} = \Uu_d; \quad \defo(\uu) = \frac{1}{2}\big(\nabla \uu + \nabla^T \uu \big);}
\end{multline}
\item
the equilibrium equations:
\begin{multline}\label{eq1:SAref}
\shoveright{\cont \in \Scb; \quad \forall \: \uu^{\ast} \in \Ucb_0, \quad \intO \Tr\big[\cont \: \defo(\uu^{\ast})\big] \dO = \intO \und{f}_d \cdot \uu^{\ast} \dO + \int_{\partial_2 \Om}\und{F}_d \cdot \uu^{\ast} \dS;}
\end{multline}
\item
the constitutive relation:
\begin{multline}\label{eq1:RDCref}
\shoveright{\cont(M) = \K \: \defo\big(\uu(M)\big) \quad \forall \: M \in \Om,}
\end{multline}
\end{itemize}

where $\Ucb = [\Hc^1(\Om)]^3$ and $\Scb = \left\{ \cont \slash \cont = \cont^{T}, \cont \in [\Lc^2(\Om)]^6 \right\}$ are the spaces of finite-energy solutions; $\Hc^1(\Om)$ stands for the standard Sobolev space of square integrable functions and first derivatives, while $\Lc^2(\Om)$ refers to the space of square integrable functions. $\Ucb_0 = \left\{ \uu^{\ast} \in \Ucb, \uu^{\ast}_{\restrictto{\partial_1 \Om}} = \und{0} \right\}$ represents the vectorial space associated to $\Ucb$, and $\K$ is the Hooke elasticity tensor.\\

In the following, the exact solution of the reference problem is designated by $(\uu, \cont)$.

\subsection{Presentation of the finite element problem}\label{2.2}

In most practical cases, the exact solution $(\uu, \cont)$ cannot be obtained explicitly. Typically, numerical methods are employed to achieve an approximation of $(\uu, \cont)$. One of the most widespread approximation methods is the finite element method (FEM) which provides a numerical solution $(\uu_h, \cont_h)$ lying in the discretized space $\Ucb_{h} \times \Scb_{h} \subset \Ucb \times \Scb$. This is defined from piecewise continuous polynomial displacement shape functions associated with a spatial discretization (finite element space mesh $\Mc_h$) of the domain $\Om$. The prescribed displacement field $\Uu_d$ is assumed to be compatible with the FE discretization. Thus, the finite element problem consists of finding a displacement/stress pair $(\uu_h (\Mu),\cont_h (\Mu)), \Mu \in \Om$, which verifies:
\begin{itemize}
\item 
the kinematic conditions:
\begin{multline}\label{eq1:CAEF}
\shoveright{\uu_h \in \Ucb_{h}; \quad {\uu_h}_{\restrictto{\partial_1 \Om}} = \Uu_d; \quad \defo(\uu_h) = \frac{1}{2}\big(\nabla \uu_h + \nabla^T \uu_h \big);}
\end{multline}
\item
the equilibrium equations:
\begin{multline}\label{eq1:SAEF}
\shoveright{\cont_h \in \Scb_h; \quad \forall \: \uu_h^{\ast} \in \Ucb_{h,0}, \quad \intO \Tr\big[\cont_h \: \defo(\uu_h^{\ast})\big] \dO = \intO \und{f}_d \cdot \uu_h^{\ast} \dO + \int_{\partial_2 \Om}\und{F}_d \cdot \uu_h^{\ast} \dS;}
\end{multline}
\item
the constitutive relation:
\begin{multline}\label{eq1:RDCEF}
\shoveright{\cont_h(M) = \K \: \defo\big(\uu_h(M)\big) \quad \forall \: M \in \Om,}
\end{multline}
\end{itemize}

where $\Ucb_{h,0} = \Ucb_h \cap \Ucb_0$.\\

In the displacement-type finite element framework, the FE solution $(\uu_h, \cont_h)$ satisfies both kinematic conditions (\ref{eq1:CAref}) and constitutive relation (\ref{eq1:RDCref}) of the reference problem, but not equilibrium equations (\ref{eq1:SAref}).\\

On the one hand, let us define the discretization error $\und{e}_h = \uu - \uu_h$, also called the exact error or true error, corresponding to the difference between the exact displacement solution and the FE one; this error enables to control and assess the numerical quality of the FE solution $(\uu_h, \cont_h)$. Usually, it is measured in terms of a suitable norm, such as the energy norm $\lnorm{\bullet}_{u, \Om} = \left( \intO \Tr\big[\K \: \defo(\bullet) \: \defo(\bullet)\big] \dO \right)^{1/2}$, which leads to a global discretization error $\lnorm{\und{e}_h}_{u, \Om}$. On the order hand, local errors $\Delta I = I(\uu) - I(\uu_h)$ can be defined if one seeks to evaluate quantities of interest, \ie functional outputs $I(\uu)$ of the solution.

\subsection{Construction of an admissible solution}\label{2.3}

In the framework of error estimation, a major interest concerns methods leading to strict and robust (global or local) error bounds for various types of mechanical problems. For this purpose, the construction of admissible solutions is required to obtain such reliable bounds \cite{Lad04,Lad83}. Furthermore, the quality of the yielded admissible solution can have a significant impact on the accuracy of the corresponding error estimate. An admissible pair, denoted $(\hat{\uu}_h, \hat{\cont}_h)$, should verify (\ref{eq1:CAref}) and (\ref{eq1:SAref}), which amounts to saying that $\hat{\uu}_h$ and $\hat{\cont}_h$ are kinematically and statically admissible, respectively. The strategy to set up an admissible solution of the reference problem relies on the post-processing of the data of both reference problem and FE solution $(\uu_h, \cont_h)$ at hand. When using a displacement-type FE method, a natural way to reconstruct $\hat{\uu}_h$ is to exploit straightforwardly the FE displacement field $\uu_h$ as this one already satisfies the kinematic conditions. Thus, one generally chooses $\hat{\uu}_h = \uu_h$ except in the case of incompressible materials, for which a special method is employed \cite{Lad92}. Conversely, the reconstruction of an admissible stress field $\hat{\cont}_h$ is more complicated as the FE stress field $\cont_h$ does not verify the equilibrium equations (it does on the FE sense only). This stage represents the key difficulty of the construction of an admissible pair. Various techniques enable to construct an admissible stress field, the best one with respect to the quality of the error estimator being based on a dual formulation of the reference problem. However, it relies on the FE calculation of two completely distinct formulations, which represents a serious drawback in terms of computing time and implementation \cite{Fra65,Kem03}. Three concurrent and well-suited techniques will be described and discussed in details in the next section.

Eventually, an appraisal of the quality of the new approximate solution $(\hat{\uu}_h, \hat{\cont}_h)$ is provided by the measure $e_{\cre}(\hat{\uu}_h, \hat{\cont}_h) = \lnorm{\hat{\cont}_h - \K \: \defo(\hat{\uu}_h)}_{\cont,\Om}$ of the error on the constitutive relation (\ref{eq1:RDCref}), with $\lnorm{\bullet}_{\cont, \Om} = \left( \intO \Tr\big[ \bullet \: \K^{-1} \: \bullet\big] \dO \right)^{1/2}$. The constitutive relation error (CRE) is linked to the global discretization error through the well-known Prager-Synge theorem \cite{Pra47} implying that the constitutive relation error is a guaranteed upper bound of the global discretization error:
\begin{align}
\lnorm{\und{e}_h}^2_{u, \Om} = \lnorm{\uu - \hat{\uu}_h}^2_{u, \Om} \leqslant \lnorm{\uu - \hat{\uu}_h}^2_{u, \Om} + \lnorm{\cont - \hat{\cont}_h}^2_{\cont,\Om} = e^2_{\cre}(\hat{\uu}_h, \hat{\cont}_h)
\end{align}
 
Therefore, the so-called constitutive relation error $e_{\cre}(\hat{\uu}_h, \hat{\cont}_h)$ represents an error estimate of the global discretization error $\lnorm{\und{e}_h}_{u, \Om}$; it can be used for goal-oriented error estimation as well \cite{Lad99a,Moi06,Cha07,Lad08,Cha08,Pan10}. Now, let us recall the basic ideas of the different techniques for constructing admissible stress fields in order to highlight the pros and cons of each method. These ideas have been more deeply presented in \cite{Lad10bis}, particularly as regards the EESPT.

\section{PRINCIPLES OF THE DIFFERENT TECHNIQUES FOR CONSTRUCTING ADMISSIBLE STRESS FIELDS}\label{3}

\subsection{Notations}\label{3.1}

Let us define $\Ec$, $\Nc$, $\Ic$ and $\Jc$ the set of elements, nodes, vertices and edges of the FE mesh $\Mc_h$, respectively. $\Ec_i^{\Nc} \subset \Ec$, $\Ec_i^{\Ic} \subset \Ec$ and $\Ec_{\Gamma}^{\Jc} \subset \Ec$ stand for the set of elements connected to node $i$, vertex $i$ and edge $\Gamma$, respectively. $\Jc_i^{\Nc} \subset \Jc$ and $\Jc_i^{\Ic} \subset \Jc$ represent the set of edges connected to node $i$ and vertex $i$, respectively. $\Nc_E^{\Ec} \subset \Nc$ and $\Nc_{\Gamma}^{\Jc} \subset \Nc$ denote the set of nodes associated with element $E$ and edge $\Gamma$, respectively. $\Ic_E^{\Ec} \subset \Ic$ and $\Ic_{\Gamma}^{\Jc} \subset \Ic$ represent the set of vertices connected to element $E$ and edge $\Gamma$. Moreover, $\und{x}_i$ designates the position of vertex $i$ in the FE mesh $\Mc_h$. Eventually, it is assumed that the FE displacement field $\und{u}_h$ belongs to $\Ucb^{p}_{h}$, where $\Ucb^{p}_{h}$ corresponds to the FE interpolation space of maximum degree $p$. $\Uc^{p}_{h}$ denotes its one-dimensional correspondent.

\subsection{The element equilibration technique (EET) - standard method}\label{3.2}

\subsubsection{Principle of the construction}\label{3.2.1}

The first technique, called the \textit{element equilibration technique} (\textit{EET}), was introduced by Ladev\`eze \cite{Lad75,Lad83}. The principle is to exploit the FE properties of the stress field $\cont_h$ through an energy relation, called \textit{prolongation condition}, which is the key point of the method. That relation links the searched admissible stress field $\hat{\cont}_h$ to the FE stress field $\cont_h$ under the form:
\begin{align}\label{eq1:prolongfort}
\intE \left(\hat{\cont}_h - \cont_h \right) \: \und{\nabla} \varphi_i \dO = \und{0} \quad \forall \: E \in \Ec, \ \forall \: i \in \Nc_E^{\Ec},
\end{align}
where $\varphi_i \in \Uc^{p}_{h}$ stands for the FE shape function associated with node $i$. This prolongation condition is a physically sound relation as it imposes that the unknown admissible stress field $\hat{\cont}_h$ provides the same work as the FE stress field $\cont_h$ for each element $E$ of the FE mesh and for all FE displacement fields. 

Classically, the EET is a quasi-explicit technique with a two-stage procedure:
\begin{enumerate}
\item[(i)]
construction of tractions $\hat{\und{F}}_h$ in equilibrium with the external loading $(\und{F}_d, \und{f}_d)$ on element boundaries $\partial E$ of the spatial mesh $\Mc_h$;\\
\item[(ii)]
calculation of an admissible stress field $\hat{\cont}_h$ in equilibrium with these equilibrated tractions $\hat{\und{F}}_h$ and body force field $\und{f}_d$ at the element level.
\end{enumerate}

\subsubsection{First stage: construction of tractions $\hat{\und{F}}_h$}\label{3.2.2}

Tractions $\hat{\und{F}}_h$ are intended to represent the stress vectors $\hat{\cont}_{h \restrictto{E}} \: \und{n}_E$ on edges $\partial E$ of element $E \in \Ec$:
\begin{align}\label{eq1:defdensites}
\hat{\cont}_{h \restrictto{E}} \: \und{n}_E = \eta_E \: \hat{\und{F}}_h \quad \textrm{on} \ \dE,
\end{align}
where $\und{n}_E$ is the unit outward normal vector to element $E$ and $\eta_E = \pm 1$ are functions ensuring the continuity of the stress vector across $\dE$.

Furthermore, these tractions $\hat{\und{F}}_h$ are built in equilibrium with the external loading $(\und{F}_d, \und{f}_d)$:
 \begin{align}
 & \eta_E \: \hat{\und{F}}_h = \und{F}_d \quad \textrm{on} \ \dE \subset \partial_2 \Om \label{eq1:eqFd} \\ 
 & \intdE \eta_E \: \hat{\und{F}}_h \cdot \uu_s \dS + \intE \und{f}_d \cdot \uu_s \dO = 0 \quad \forall \: \uu_s \in \Ucb_{S \restrictto{E}}, \label{eq1:eqfd}
 \end{align}
where $\Ucb_{S \restrictto{E}}$ is the set of rigid body displacement fields over element $E$.

The procedure for calculating tractions $\hat{\und{F}}_h$ on the element boundaries of the spatial mesh $\Mc_h$ is quasi-explicit. First, tractions $\hat{\und{F}}_h$ are set equal to the traction force density $\und{F}_d$ on the subdomain $\partial_2 \Om$: $\eta_E \: \hat{\und{F}}_h = \und{F}_d$ on $\dE \subset \partial_2 \Om$. Second, by applying integration by parts (Green's theorem) before using definition (\ref{eq1:defdensites}) of tractions $\hat{\und{F}}_h$ as well as the strong form of the problem associated to equilibrium equations (\ref{eq1:SAref}) satisfied by $\hat{\cont}_h$, the prolongation condition (\ref{eq1:prolongfort}) can be recasted in the following form:
\begin{align}\label{eq1:densites}
\intdE \eta_E \: \hat{\und{F}}_h \: \varphi_i \dS = \intE \left(\cont_h \: \und{\nabla} \varphi_i - \und{f}_d \: \varphi_i \right) \dO = \und{\Qc}_E^{(i)} \quad \forall \: E \in \Ec, \ \forall \: i \in \Nc_E^{\Ec},
\end{align}

\begin{remark}
The derivation of (\ref{eq1:densites}) also assumes some regularity for $\hat{\cont}_{h \restrictto{E}}$ over each element $E \in \Ec$: either $\hat{\cont}_{h \restrictto{E}}$ is continuous inside $E$ or, if a discontinuity occurs along a path $\Gamma \in E$, the normal component $\hat{\cont}_{h \restrictto{\Gamma}} \: \und{n}_{\Gamma}$ along $\Gamma$ is continuous.
\end{remark}

It is worthy noticing that (\ref{eq1:densites}) enforces the equilibrium condition (\ref{eq1:eqfd}) of tractions $\hat{\und{F}}_h$ with the body force field $\und{f}_d$ over each element $E \in \Ec$, as (\ref{eq1:densites}) is equivalent to
\begin{align}\label{eq2:densites}
\intdE \eta_E \: \hat{\und{F}}_h \cdot \uu_h^{\ast} \dS = \intE \left(\Tr\big[\cont_h \: \defo(\uu_h^{\ast})\big] - \und{f}_d \cdot \uu_h^{\ast} \right) \dO \quad \forall \: \uu_h^{\ast} \in \Ucb^p_{h \restrictto{E}},
\end{align}
and $\Ucb_{S \restrictto{E}} \subset \Ucb^{p}_{h \restrictto{E}}$.

Thus, the prolongation condition (\ref{eq1:densites}) defines local problems $\Pc_i^{\Nc}$ associated with each node $i \in \Nc$. Problem $\Pc_i^{\Nc}$ associated with node $i \in \Nc$ is a linear system that reads:
\begin{align}\label{eq1:systproj}
\displaystyle\sum_{\Gamma \in \Jc_i^{\Nc} \cap \dE} \eta_E \: \hat{\und{b}}_{\restrictto{\Gamma}}^{(i)} = \und{\Qc}_E^{(i)} \quad \forall \: E \in \Ec_i^{\Nc},
\end{align}
where the unknown quantity $\displaystyle\hat{\und{b}}_{\restrictto{\Gamma}}^{(i)}$ over the edge $\Gamma \in \Jc_i^{\Nc}$ is the projection of traction $\hat{\und{F}}_{h \restrictto{\Gamma}}$ over the FE shape function $\varphi_i$:
\begin{align}\label{eq1:defproj}
\displaystyle\hat{\und{b}}_{\restrictto{\Gamma}}^{(i)} = \int_{\Gamma} \hat{\und{F}}_h \: \varphi_i \dS.
\end{align}

Solvability of local problems $\Pc_i^{\Nc}$ (\ref{eq1:systproj}) depends on the type of node $i$ considered. Indeed, existence of a solution for system (\ref{eq1:systproj}) most of the time requires the verification of a compatibility condition, that reads:
\begin{align}
& \sum_{E \in \Ec_i^{\Nc}} \und{\Qc}_E^{(i)} = \und{0} \qquad \textrm{for an internal node} \ i \in \Om, \label{eq1:compatibilite} \\
& \sum_{E \in \Ec_i^{\Nc}} \und{\Qc}_E^{(i)} = \int_{\Gamma \in \Jc_i^{\Nc} \cap \partial_2 \Om} \und{F}_d \: \varphi_i \dS \qquad \textrm{for a node} \ i \subset \partial_2 \Om. \label{eq2:compatibilite}
\end{align}

Noticing that: 
\begin{align}\label{eq1:sumproj}
\sum_{E \in \Ec_i^{\Nc}} \und{\Qc}_E^{(i)} = \intO \left(\cont_h \: \und{\nabla} \varphi_i - \und{f}_d \: \varphi_i \right) \dO,
\end{align}
the compatibility conditions (\ref{eq1:compatibilite}) and (\ref{eq2:compatibilite}) are ensured by the FE equilibrium equations (\ref{eq1:SAEF}) satisfied by the FE stress field $\cont_h$, as the vectorial equation (\ref{eq1:sumproj}) is linked after some straightforward computations to the scalar equilibrium equations (\ref{eq1:SAEF}).\\

If local problems $\Pc_i^{\Nc}$ (\ref{eq1:systproj}) can admit several solutions, the uniqueness of the solution is ensured by performing the least-squares minimization of a cost function of the form \cite{Lad04}:
\begin{align}\label{eq1:fonctioncoutEET}
J(\hat{\und{b}}^{(i)}) = \frac{1}{2} \sum_{\Gamma \in \Jc_i^{\Nc}} ( \hat{\und{b}}^{(i)} - \und{b}^{(i)} )_{\restrictto{\Gamma}}^2,
\end{align}
where the known quantity $\displaystyle\und{b}_{\restrictto{\Gamma}}^{(i)}$ represents the projection of the FE stress vector $\cont_h \: \und{n}$, over the edge $\Gamma \in \Jc_i^{\Nc}$, on the FE shape function $\varphi_i$:
\begin{itemize}
\item
for an internal edge $\Gamma$ connected to elements $E$ and $E'$,
\begin{align*}
\displaystyle\und{b}^{(i)}_{\restrictto{\Gamma}} = \frac{1}{2} \int_{\Gamma} (\frac{1}{\eta_{E}} \cont_{h \restrictto{E}} \: \und{n}_E + \frac{1}{\eta_{E'}} \cont_{h \restrictto{E'}} \: \und{n}_{E'}) \: \varphi_i \dS;
\end{align*}
\item
for an edge $\Gamma \subset \partial_1 \Om$ connected to element $E$,
\begin{align*}
\displaystyle\und{b}^{(i)}_{\restrictto{\Gamma}} = \int_{\Gamma} \frac{1}{\eta_{E}} \cont_{h \restrictto{E}} \: \und{n}_E \: \varphi_i \dS;
\end{align*}
\item
for an edge $\Gamma \subset \partial_2 \Om$ connected to element $E$,
\begin{align*}
\displaystyle\und{b}^{(i)}_{\restrictto{\Gamma}} = \int_{\Gamma} \frac{1}{\eta_{E}} \und{F}_d \: \varphi_i \dS.
\end{align*}
\end{itemize}
Other cost functions will be considered in \Sect{5} to study the influence of its choice on the quality of the error estimate.\\

Eventually, densities $\hat{\und{F}}_h$ along edges $\Gamma \in \Jc$ are recovered from $\displaystyle\hat{\und{b}}^{(i)}$, using the same interpolation degree as the FE displacement field $\uu_h$. In other words, one seeks $\displaystyle\hat{\und{F}}_{h \restrictto{\Gamma}} \in \Ucb^{p}_{h \restrictto{\Gamma}}$ under the form: 
\begin{align}\label{eq1:interpol_Fh_hat}
\hat{\und{F}}_{h \restrictto{\Gamma}} = \sum_{j \in \Nc_{\Gamma}^{\Jc}} \hat{\und{F}}_j^{\Gamma} \: \varphi_{j \restrictto{\Gamma}}.
\end{align}
Therefore, the evaluation of tractions $\displaystyle\hat{\und{F}}_{h \restrictto{\Gamma}}$ along each edge $\Gamma \in \Jc$ requires the solution of a set of linear local problems at the edge level involving the projections $\displaystyle\hat{\und{b}}_{\restrictto{\Gamma}}^{(i)}$ associated to each node $i \in \Nc_{\Gamma}^{\Jc}$.

\subsubsection{Second stage: calculation of an admissible stress field $\hat{\cont}_h$ at the element scale}\label{3.2.3}

Once the flux-equilibration procedure has been performed, one seeks the local restriction $\hat{\cont}_{h \restrictto{E}}$ of an admissible stress field $\hat{\cont}_h$ to each element $E \in \Ec$ as the solution of an equilibrium local problem $\Pc_E^{\Ec}$ where the previously calculated tractions $\displaystyle\hat{\und{F}}_h$ play the role of the external loading over the element boundaries:
\begin{equation}\label{eq1:pblocal}
\hat{\cont}_{h \restrictto{E}} \in \Scb_{\hat{\und{F}}_h} \iff \left\{
\begin{aligned}
& \hat{\cont}_{h \restrictto{E}} \in \Scb \\
& \und{\diver} \: \hat{\cont}_{h \restrictto{E}} + \und{f}_d = \und{0} \quad \textrm{in} \ E \\
& \hat{\cont}_{h \restrictto{E}} \: \und{n}_E = \eta_E \: \hat{\und{F}}_h \quad \textrm{on} \ \dE
\end{aligned}
\right.
\end{equation}

It is important to recall that the construction of tractions $\displaystyle\hat{\und{F}}_h$ leads to equilibrated local problems $\Pc_E^{\Ec}$ (\ref{eq1:pblocal}) ensuring their solvability, since $\hat{\cont}_{h \restrictto{E}}$ is assumed to have no discontinuity inside each element $E \in \Ec$.
In practice, a dual formulation (displacement-type FEM) of local static problems $\Pc_E^{\Ec}$ (\ref{eq1:pblocal}) combined with a $p$-refinement or $h$-refinement of the FEM is used to obtain a numerical approximation of $\hat{\cont}_{h \restrictto{E}}$. In other words, the discretization of each element $E \in \Ec$ consists of considering either a single element $E$ along with an interpolation of degree $p+k$, where $p$ denotes the degree of the FE analysis and $k$ an additional degree, or a subdivision of element $E$ along with an interpolation of degree $p$. The $p$-refinement technique, introduced by Babu\u{s}ka and Strouboulis \cite{Bab94}, has a higher convergence rate than the $h$-refinement one and circumvents the need to generate a new refined mesh. Using an extra-degree $k=3$ yields a very good and accurate approximation of a statically admissible stress field, thus leading to an error estimate of good quality \cite{Coo99}. The $p$-refinement technique could also be combined with a $h$-refinement technique. The interested reader can refer to \cite{Lad04} for more information.

\begin{remark}
The stress field $\hat{\cont}_{h \restrictto{E}}$ obtained using a $p$-refinement technique is not strictly admissible in the sense that local problems $\Pc_E^{\Ec}$ (\ref{eq1:pblocal}) are not solved exactly; such a stress field leads to only asymptotic error bounds with respect to a refined solution.
\end{remark}

\begin{remark}
Local problems $\Pc_E^{\Ec}$ (\ref{eq1:pblocal}) can however be solved exactly: a first procedure consists of constructing an admissible stress field analytically; this construction requires a subsplitting of elements, the statically admissible stress field being searched using a piecewise polynomial basis \cite{Lad04,Lad96}. This analytical construction yields strictly admissible stress fields provided that body force field $\und{f}_d$ is a polynomial with a degree compatible with that of $\hat{\cont}_{h \restrictto{E}}$. Another procedure consists of using a stress-type FEM. Both procedures yield admissible stress fields and therefore strict error bounds with respect to the exact solution.
\end{remark}

The EET has a very attractive feature, namely its affordable computational cost, as local problems $\Pc_i^{\Nc}$ (\ref{eq1:systproj}) and $\Pc_E^{\Ec}$ (\ref{eq1:pblocal}) are both reasonably costly. Despite this advantage, the flux-equilibration procedure needed to properly set boundary conditions for local problems $\Pc_E^{\Ec}$ (\ref{eq1:pblocal}) is complex to implement, which strongly limits its practical use into existing FE codes. Nevertheless, it has already been implemented into the industrial software SAMCEF \cite{Samtech}.

\subsection{The star-patch equilibration technique (SPET) - flux-free method}\label{3.3}

The second method, called the \textit{star-patch equilibration technique} (\textit{SPET}), was formerly developed for fluid mechanics \cite{Car00,Mor03,Pru04} and has been revisited and adapted to solid mechanics by Par\'es, D\'iez and Huerta \cite{Par06} under the name \textit{flux-free technique}. This technique does not require any flux recovery or flux splitting technique as there is no need to perform flux equilibration; this constitutes a serious advantage for implementation. The novelty of this technique is the introduction of a partition of unity which allows to boil down to self-equilibrated local problems defined on sets of elements, also called \textit{patches} or \textit{stars}.

\subsubsection{Definition of self-equilibrated local problems over patches of elements}\label{3.3.1}

The global problem defining the discretization error $\und{e}_h = \uu - \uu_h$ can be obtained from the equilibrium equations (\ref{eq1:SAref}) replacing $\uu$ by $\und{e}_h + \uu_h$. It reads:

Find $\und{e}_h \in \Ucb_0$ such that:\\
\begin{equation}\label{eq1:residus}
\begin{aligned}
 \intO \Tr\big[\K \: \defo(\und{e}_h) \: \defo(\uu^{\ast})\big] \dO & = - \intO \Tr\big[\K \: \defo(\uu_h) \: \defo(\uu^{\ast})\big] \dO + \intO \und{f}_d \cdot \uu^{\ast} \dO + \int_{\partial_2 \Om}\und{F}_d \cdot \uu^{\ast} \dS \\
 & = \Rc_h(\uu^{\ast}) \quad \forall \: \uu^{\ast} \in \Ucb_0,
\end{aligned}
\end{equation}
where $\Rc_h$ stands for the weak residual functional associated with the FE approximation $\uu_h$. This residual equation conveys the failure to comply with the FE equilibrium equations (\ref{eq1:SAEF}).\\

The main idea consists of introducing at this stage the partition of unity defined by the linear FE shape functions $\la_i \in \Uc_h^1$ based on vertices $i \in \Ic$, which are such that:
\begin{align*}
\la_i (\und{x}_j) = \delta_{ij} \quad \forall \: j \in \Ic; \quad \sum_{i \in \Ic} \la_i(\und{x}) = 1 \quad \textrm{for} \ \und{x} \in \Om.
\end{align*}
Global problem (\ref{eq1:residus}) can thus be reformulated as:
\begin{align}\label{eq2:residus}
 \intO \Tr\big[\K \: \defo(\und{e}_h) \: \defo(\uu^{\ast})\big] \dO = \sum_{i \in \Ic} \Rc_h( \la_i \: \uu^{\ast} ) \quad \forall \: \uu^{\ast} \in \Ucb_0.
\end{align}

\begin{figure}
\centering\includegraphics[scale=0.25]{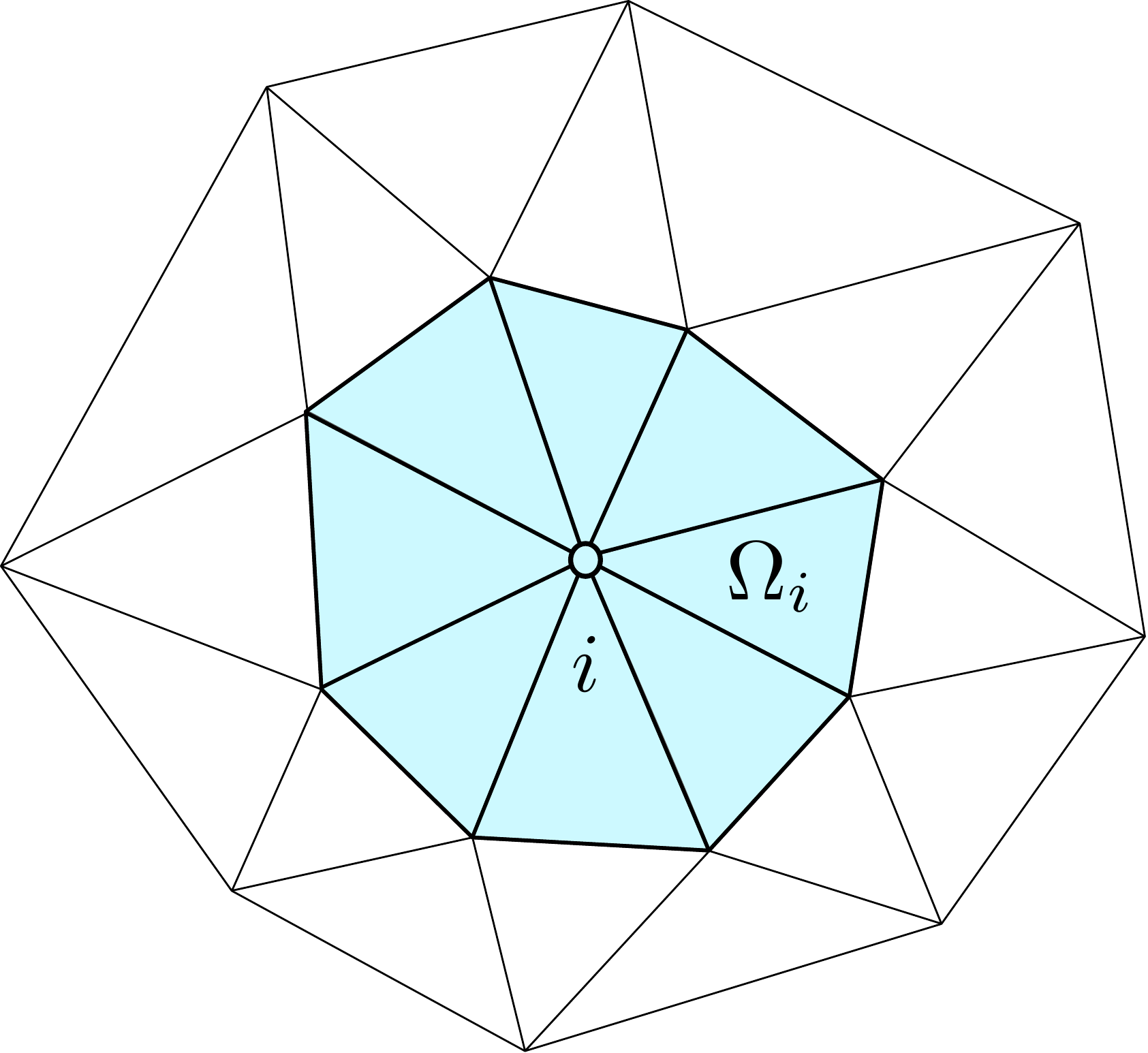}
\caption{Patch $\Om_i$ associated with vertex $i$.}\label{fig1:patch}
\end{figure}

Let us define the patch $\Om_i$ as the support of $\la_i$; in other words, $\Om_i$ designates the support of the set $\Ec_i^{\Ic}$ of elements connected to vertex $i \in \Ic$ (see \Fig{fig1:patch}).

Therefore, the SPET consists of solving a set of local problems in each patch $\Om_i$:

Find $\und{e}_i \in \Ucb_{0 \restrictto{\Om_i}}$ such that:
\begin{align}\label{eq1:pbpatch}
 \intOi \Tr\big[\K \: \defo(\und{e}_i) \: \defo(\uu^{\ast})\big] \dO = \Rc_h( \la_i \: \uu^{\ast} ) \quad \forall \: \uu^{\ast} \in \Ucb_{0 \restrictto{\Om_i}},
\end{align}
where $\Ucb_{0 \restrictto{\Om_i}}$ is the restriction of $\Ucb_0$ to patch $\Om_i$: $\Ucb_{0 \restrictto{\Om_i}} = \left\{ \uu^{\ast} \in \Ucb_{\restrictto{\Om_i}}, \uu^{\ast}_{\restrictto{\partial_1 \Om \cap \Om_i}} = \und{0} \right\}$.\\

Now, let us recall the key property of the Galerkin approach, which is the Galerkin orthogonality:
\begin{align}\label{eq1:orthogalerkin}
\Rc_h(\uu^{\ast}_h) = 0 \quad \forall \: \uu^{\ast}_h \in \Ucb_{h, 0}^p,
\end{align}
where $\Ucb_{h, 0}^p = \Ucb_{h}^p \cap \Ucb_0$.

Let us note that the weak residual $\Rc_h(\la_i \: \bullet)$ does not allow to ensure the solvability of local problems (\ref{eq1:pbpatch}). Indeed, existence of a solution of (\ref{eq1:pbpatch}) is guaranteed if and only if the following compatibility condition holds: $\Rc_h(\la_i \: \uu_s) = 0 \quad \forall \: \uu_s \in \Ucb_{S,0 \restrictto{\Om_i}}$, where $\Ucb_{S,0 \restrictto{\Om_i}}$ is the restriction of $\Ucb_{S,0}$ to patch $\Om_i$, \ie the set of rigid body displacement fields defined over $\Om_i$ which vanish on $\partial_1 \Om \cap \Om_i$. However, this condition is not verified for an interpolation degree $p = 1$, as $\la_i \: \uu_s \in \Ucb_{h, 0}^2 \not \subset \Ucb_{h, 0}^1$.\\

In \cite{Par06}, the projector $\Pi \colon \Ucb_{0} \to \Ucb_{h, 0}^1$ whose image is in the space of linear FE shape functions was defined. By using Galerkin orthogonality (\ref{eq1:orthogalerkin}), the right-hand side term of (\ref{eq1:pbpatch}) can be replaced by $\Rc_h( \la_i \: (\uu^{\ast} - \Pi \uu^{\ast}) )$ leading to a new set of local problems $\Pc_i^{\Ic}$ over each patch $\Om_i$ in the form:

Find $\und{e}_i \in \Ucb_{0 \restrictto{\Om_i}}$ such that:
\begin{align}\label{eq2:pbpatch}
 \intOi \Tr\big[\K \: \defo(\und{e}_i) \: \defo(\uu^{\ast})\big] \dO = \Rc_h( \la_i \: (\uu^{\ast} - \Pi \uu^{\ast}) ) \quad \forall \: \uu^{\ast} \in \Ucb_{0 \restrictto{\Om_i}}.
\end{align}

By observing that:
\begin{align}\label{eq1:eqR_h}
\Pi \uu_s = \uu_s \quad \forall \: \uu_s \in \Ucb_{S,0 \restrictto{\Om_i}} \implies \Rc_h( \la_i \: (\uu_s - \Pi \uu_s) ) = 0 \quad \forall \: \uu_s \in \Ucb_{S,0 \restrictto{\Om_i}},
\end{align}
local problems $\Pc_i^{\Ic}$ (\ref{eq2:pbpatch}) are well-posed and solvable, the projector being the key ingredient to ensure their solvability, especially for an interpolation degree $p = 1$. Let us note that for scalar (\ie one-dimensional) problems and for two- or three-dimensional mechanical problems with high-order interpolation degree (at least quadratic), it is not necessary to introduce the projector $\Pi$ into the right-hand side term of (\ref{eq1:pbpatch}) to achieve the solvability of local problems $\Pc_i^{\Ic}$ (\ref{eq2:pbpatch}).\\

Therefore, local problems (\ref{eq2:pbpatch}) defined at the patch scale are self-equilibrated so that their solution does not require the construction of tractions or flux jumps along element boundaries and there is no need to perform any flux equilibration to achieve equilibrium, which is a very interesting feature for implementation purposes. In practice, the solution of local problems (\ref{eq2:pbpatch}) is computed using a $p$-refinement or $h$-refinement over patches $\Om_i$. The calculation is classically performed using the original FE mesh $\Mc_h$ with a $p + 3$ discretization over each patch $\Om_i$ in order to obtain a fairly good approximation of the solution of (\ref{eq2:pbpatch}).

\subsubsection{Construction of both an admissible stress field and a global estimate}\label{3.3.2}

Summing (\ref{eq2:pbpatch}) for all vertices $i \in \Ic$ leads to a relation between the discretization error $\und{e}_h$ and numerical solutions $\und{e}_i$:
\begin{align}
\sum_{i \in \Ic} \intOi \Tr\big[\K \: \defo(\und{e}_i) \: \defo(\uu^{\ast})\big] \dO = \sum_{i \in \Ic} \Rc_h(\la_i \: (\uu^{\ast} - \Pi \uu^{\ast}) ) = \intO \Tr\big[\K \: & \defo(\und{e}_h) \: \defo(\uu^{\ast})\big] \dO
\end{align}

Consequently, numerical solutions $\und{e}_i$ obtained from (\ref{eq2:pbpatch}) allow to define a global error estimate, which is an overestimation of the energy norm of the discretization error:
\begin{align}
\lnorm{\und{e}_h}_{u, \Om} \: \leqslant \: \left( \sum_{E\in \Ec} \lnorm{\sum_{i \in \Ic_E^{\Ec}} \und{e}_{i \restrictto{E}}}^2_{u, E} \right)^{1/2},
\end{align}
with $\displaystyle\lnorm{\bullet}_{u, E} = \left( \intE \Tr\big[\K \: \defo(\bullet_{\restrictto{E}}) \: \defo(\bullet_{\restrictto{E}})\big] \dO \right)^{1/2} $.

Equivalently, one can define an admissible stress field over each element $E \in \Ec$:
\begin{align}
 \hat{\cont}_{h \restrictto{E}} = \cont_{h \restrictto{E}} + \K \: \defo \left(\sum_{i \in \Ic_E^{\Ec}} \und{e}_{i \restrictto{E}}\right).
\end{align} 

Finally, the SPET evades the main disadvantage of the EET, because it does not require the construction of equilibrated tractions making its implementation less cumbersome. On the other hand, local problems $\Pc_i^{\Ic}$ (\ref{eq2:pbpatch}) are defined at the patch scale and, therefore, involve a large number of degrees of freedom, particularly in 3D, compared with local problems $\Pc_E^{\Ec}$ (\ref{eq1:pblocal}) defined at the element scale for the EET. Thus, the calculation of an admissible stress field using the SPET may result in simply defined but costly computations.

\begin{remark}
Local problems $\Pc_i^{\Ic}$ (\ref{eq2:pbpatch}) defined on patches can be solved using a dual approach, that is a stress-type FEM \cite{Cot09}.
\end{remark}

\subsection{The element equilibration and star-patch technique (EESPT) - hybrid method}\label{3.4}

\subsubsection{Principle of the construction}\label{3.4.1}

This new technique, based on recent works of Ladev\`eze \textit{et al} \cite{Lad10bis}, is a hybrid method insofar as it takes advantage of the ingredients of both EET and SPET methods. As for the EET, the procedure to construct an admissible stress field involves two main stages:
\begin{enumerate}
\item[(i)]
construction of tractions $\hat{\und{F}}_h$ in equilibrium with the external loading $(\und{F}_d, \und{f}_d)$ on element boundaries $\partial E$ of the spatial mesh $\Mc_h$;\\
\item[(ii)]
construction of an admissible stress field $\hat{\cont}_h$ solution of the static local problem $\Pc_E^{\Ec}$ (\ref{eq1:pblocal}) over each element $E \in \Ec$.
\end{enumerate}
The second stage is similar to that of the EET. As previously explained, those local problems $\Pc_E^{\Ec}$ are solved numerically by using a dual formulation of (\ref{eq1:pblocal}) and the FEM through a refined mesh ($p$-refinement); one can refer to \Sect{3.2.3} for details.
The main thrust relates to the first stage, which aims at constructing equilibrated tractions, in order to facilitate its practical implementation. The principle is to exploit not only the prolongation condition used in the EET, but also the partition of unity involved in the SPET.

\subsubsection{Definition of a set of local problems over patches of elements}\label{3.4.2}

Prolongation condition (\ref{eq1:prolongfort}), in which the FE shape functions $\varphi_i$ belong to $\Uc^{p}_{h}$, can be reformulated in the global form:
\begin{align}\label{eq2:prolongfort}
\intO \Tr\big[\left( \hat{\cont}_h - \cont_h \right) \: \defo(\und{v}_h^{\ast}) \big] \dO = \sum_{E \in \Ec} \intE \Tr\big[\left( \hat{\cont}_h - \cont_h \right) \: \defo(\und{v}_h^{\ast}) \big] \dO = 0 \quad \forall \: \und{v}_h^{\ast} \in \Vcb_h^p,
\end{align}
where $\Vcb_h^p$ designates the space of polynomial functions which are continuous over each element $E \in \Ec$ but possibly discontinuous across the inter-element edges.

One can replace the vectorial space $\Vcb_h^p$ of degree $p$ by the vectorial space $\Vcb_h^q$ of degree $q$ such that $1 \leqslant q \leqslant p$, since this one suffices to ensure the equilibrium condition (\ref{eq1:eqfd}). Let us note that using $q < p$ is not equivalent to enforcing prolongation condition (\ref{eq1:prolongfort}). For the sake of simplicity and affordable computational cost, one limits to the case $q = 1$, but the method can be generalized for functions $\und{v}_h^{\ast} \in \Vcb_h^q$ of degree $q$ with $1 \leqslant q \leqslant p$.

By considering linear FE shape functions $\varphi_i \in \Uc_{h}^1$ only, (\ref{eq2:prolongfort}) becomes:
\begin{align}\label{eq3:prolongfort}
\intO \Tr\big[\left( \hat{\cont}_h - \cont_h \right) \: \defo(\und{v}_h^{\ast}) \big] \dO = \sum_{E \in \Ec} \intE \Tr\big[\left( \hat{\cont}_h - \cont_h \right) \: \defo(\und{v}_h^{\ast}) \big] \dO = 0 \quad \forall \: \und{v}_h^{\ast} \in \Vcb_h^1.
\end{align}

Using the weak equilibrium between $\hat{\cont}_h$, $\hat{\und{F}}_h$ and $\und{f}_d$ gives:
\begin{align}\label{eq3:densites}
\sum_{E \in \Ec} \left[ \intdE \eta_E \: \hat{\und{F}}_h \cdot \und{v}_h^{\ast} \dS - \int_E \left( \Tr\left[\cont_h \: \defo(\und{v}_h^{\ast})\right] - \und{f}_d \cdot \und{v}_h^{\ast} \right) \dO \right] = 0 \quad \forall \: \und{v}_h^{\ast} \in \Vcb_h^1.
\end{align}

Then, incorporating the partition of unity defined by the linear FE shape functions $\la_i \in \Uc_h^1$ into (\ref{eq3:densites}) leads to:
\begin{align}\label{eq4:densites}
\sum_{E \in \Ec_i^{\Ic}} \left[ \intdE \eta_{E} \: \la_i \: \hat{\und{F}}_h \cdot \und{v}_h^{\ast} \dS - \int_E \left( \Tr\left[\cont_h \: \defo(\la_i \: \und{v}_h^{\ast})\right] - \und{f}_d \cdot \la_i \: \und{v}_h^{\ast} \right) \dO \right] = 0 \quad \forall \: \und{v}_h^{\ast} \in \Vcb_h^1.
\end{align}

In order to offer better flexibility, we will consider the following set of local problems $\Pc_i^{\Ic}$ introduced in \cite{Lad10bis} and defined over each patch $\Om_i$ by observing that $\Om_i$ is the support of $\la_i$:

Find $\la_i \: \hat{\und{F}}_h^{(i)}$ such that:
\begin{align}\label{eq1:pbdensitespatch}
\sum_{E \in \Ec_i^{\Ic}} \left[ \intdE \eta_{E} \: \la_i \: \hat{\und{F}}_h^{(i)} \cdot \und{v}_h^{\ast} \dS \right] = \Qc_{\Om_i}(\la_i \: \und{v}_h^{\ast}) \quad \forall \: \und{v}_h^{\ast} \in \Vcb_h^1,
\end{align}
where $\displaystyle \Qc_{\Om_i}(\la_i \: \und{v}_h^{\ast}) = \intOi \left( \Tr\left[\cont_h \: \defo(\la_i \: \und{v}_h^{\ast})\right] - \und{f}_d \cdot \la_i \: \und{v}_h^{\ast} \right) \dO$.

Let us note that the sum of local problems (\ref{eq1:pbdensitespatch}) (or (\ref{eq4:densites})) for all vertices $i \in \Ic$ leads to equation (\ref{eq3:densites}), which shows that global equilibrium is ensured.

By noticing that solutions $\la_i \: \hat{\und{F}}_h^{(i)}$ of (\ref{eq1:pbdensitespatch}) are nonzero only along the edges connected to vertex $i$, \ie $\Gamma \in \Jc_i^{\Ic}$, local problems $\Pc_i^{\Ic}$ (\ref{eq1:pbdensitespatch}) can be rewritten as:

Find $\la_i \: \hat{\und{F}}_h^{(i)}$ such that:
\begin{align}\label{eq2:pbdensitespatch}
 \sum_{\Gamma \in \Jc_i^{\Ic}} \int_{\Gamma} \la_i \: \hat{\und{F}}_h^{(i)} \cdot \left( \sum_{E \in \Ec_{\Gamma}^{\Jc}} \eta_{E} \: {\und{v}_h^{\ast}}_{\restrictto{E}} \right) \dS = \Qc_{\Om_i}(\la_i \: \und{v}_h^{\ast}) \quad \forall \: \und{v}_h^{\ast} \in \Vcb_h^1.
\end{align}

\begin{figure}
\centering\includegraphics[scale = 0.4]{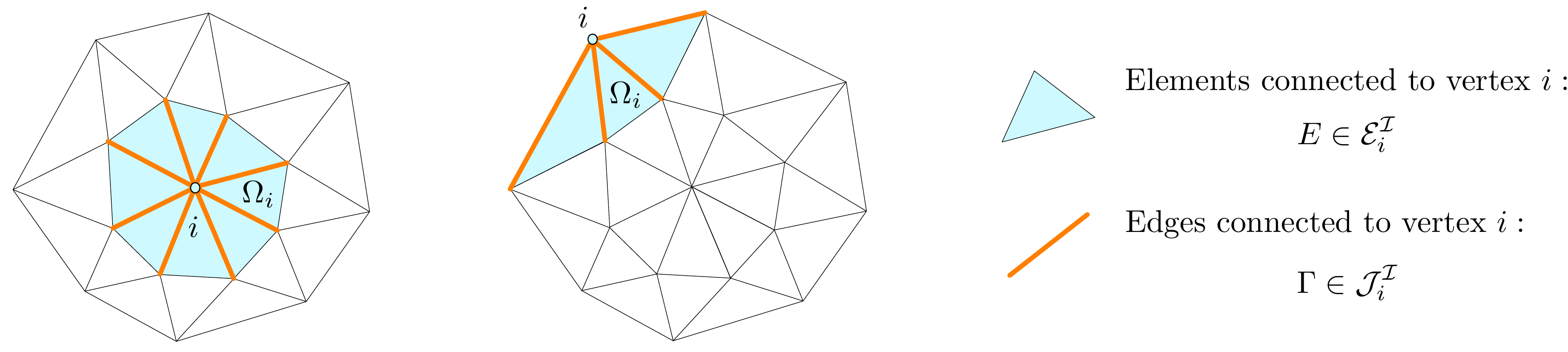}
\caption{Elements and edges connected to vertex $i$.}\label{fig1:patchdouble_black_with_legend}
\end{figure}

Both sets $\Jc_i^{\Ic}$ and $\Ec_{i}^{\Ic}$ are represented in \Fig{fig1:patchdouble_black_with_legend}.

Consequently, once local problems (\ref{eq2:pbdensitespatch}) have been solved, tractions $\hat{\und{F}}_h$ are recovered on each edge $\Gamma \in \Jc$ by summing the solutions $\la_i \: \hat{\und{F}}_h^{(i)}$ of local problems (\ref{eq2:pbdensitespatch}) for all the vertices connected to the edge $\Gamma$, \ie for $i \in \Ic_{\Gamma}^{\Jc}$:
\begin{align}\label{eq1:interpolationFh}
\hat{\und{F}}_{h \restrictto{\Gamma}} = \sum_{i \in \Ic_{\Gamma}^{\Jc}} ( \la_i \: \hat{\und{F}}_h^{(i)} )_{\restrictto{\Gamma}}.
\end{align}

In practice, local quantities $\la_i \: \hat{\und{F}}_h^{(i)}$ are sought along each edge $\Gamma \in \Jc$ in the space of piecewise continuous polynomial functions of degree $p$ defined on the element edges. Indeed, more precisely, $\displaystyle \la_i \: \hat{\und{F}}_h^{(i)} \in \Ucb^{p}_{h \restrictto{\Gamma}}$ in order to recover $\hat{\und{F}}_h$ with the same interpolation degree as the FE displacement field $\uu_h \in \Ucb_h^p$.

Furthermore, let us recall that conditions $\eta_E \: \hat{\und{F}}_h = \und{F}_d$ have to be enforced over edges $\Gamma \subset \partial_2 \Om$.

\subsubsection{Solvability and uniqueness of solution $\la_i \: \hat{\und{F}}_h^{(i)}$ for local problem $\Pc_i^{\Ic}$}\label{3.4.3}

First, let us address the question of solvability of local problems $\Pc_i^{\Ic}$. A detailed analysis of local problems $\Pc_i^{\Ic}$ over $\Om_i$ is carried out in \cite{Lad10bis}. Let us recall the main points. By considering the set ${\bar{\Ucb}_{h,0 \restrictto{\Om_i}}^1}$ of polynomial functions $\und{v}_h^{\ast} \in \Vcb_h^1$ which are continuous across edges $\Gamma \in \Jc_i^{\Ic}$ and which vanish on $\Gamma \in \Jc_i^{\Ic} \cap \dOm$, \ie such that $\und{v}_h^{\ast} = \und{0}$ along $\Gamma \in \Jc_i^{\Ic} \cap \dOm$, the existence of a solution of problem (\ref{eq2:pbdensitespatch}) is ensured for $p \geqslant 2$. However, in the case $p = 1$, local problem $\Pc_i^{\Ic}$ over $\Om_i$ is replaced by:

Find $\la_i \: \hat{\und{F}}_h^{(i)}$ such that:
\begin{align}\label{eq3:pbdensitespatch}
\sum_{\Gamma \in \Jc_i^{\Ic}} \int_{\Gamma} \la_i \: \hat{\und{F}}_h^{(i)} \cdot \left( \sum_{E \in \Ec_{\Gamma}^{\Jc}} \eta_{E} \: {\und{v}_h^{\ast}}_{\restrictto{E}} \right) \dS = \Qc_{\Om_i}(\la_i \: (\und{v}_{h}^{\ast} - \Pi^i \und{v}_{h}^{\ast})) \quad \forall \: \und{v}_h^{\ast} \in \Vcb_h^1,
\end{align}
where $\displaystyle \Qc_{\Om_i}(\la_i \: (\und{v}_{h}^{\ast} - \Pi^i \und{v}_{h}^{\ast})) = \Qc_{\Om_i}(\la_i \: \und{v}_{h}^{\ast}(\und{x}_i)) = \intOi \left( \Tr\left[\cont_h \: \defo(\la_i \: \und{v}_h^{\ast}(\und{x}_i))\right] - \und{f}_d \cdot \la_i \: \und{v}_h^{\ast}(\und{x}_i) \right) \dO$.
Operator $\Pi^i \colon \Vcb_h^1 \to \Vcb_{h,i}^1$ designates the projector onto the space $\Vcb_{h,i}^1$ of functions $\und{v}_h^{\ast} \in \Vcb_h^1$ defined over $\Om_i$ which vanish at vertex $i$. Mathematically, it reads $\Pi^i \und{v}_{h \restrictto{E}}^{\ast} = \und{v}_{h \restrictto{E}}^{\ast} - \und{v}_{h \restrictto{E}}^{\ast}(\und{x}_i) \quad \forall \: E \in \Ec_i^{\Ic}, \ \forall \: \und{v}_{h}^{\ast} \in \Vcb_h^1$. The solvability of local problem $\Pc_i^{\Ic}$ defined over patch $\Om_i$ by (\ref{eq3:pbdensitespatch}) is ensured and those local problems introduced for $p = 1$ implicitly require that (see \cite{Lad10bis}):
\begin{align}\label{eq3:wh}
\sum_{\Gamma \in \Jc_i^{\Ic}} \int_{\Gamma} \la_i \: \hat{\und{F}}_h^{(i)} \: \la_j \dS = \und{0} \quad \forall \: j \in \Ic_{E}^{\Ec} \: \textrm{and} \: j \neq i.
\end{align}

\begin{remark}
In the case $q = p$, $\und{v}_h^{\ast} \in \Vcb_h^p$ and local problems (\ref{eq2:pbdensitespatch}) are not solvable. Therefore, the modification introduced in (\ref{eq3:pbdensitespatch}) should be used to ensure solvability.
\end{remark}

Now, let us focus on the uniqueness of solution $\la_i \: \hat{\und{F}}_h^{(i)}$ of (\ref{eq2:pbdensitespatch}) (or (\ref{eq3:pbdensitespatch})). This one is obtained by performing a least-squares minimization step over each patch $\Om_i$ involving a cost function whose form is similar to (\ref{eq1:fonctioncoutEET}):
\begin{align}\label{eq1:fonctioncoutEESPT}
J_{\Om_i}(\la_i \: \hat{\und{F}}_h^{(i)}) = \frac{1}{2} \sum_{\Gamma \in \Jc_i^{\Ic}} ( \la_i \: \hat{\und{F}}_h^{(i)} - \la_i \: \und{F}_h^{(i)} )_{\restrictto{\Gamma}}^2.
\end{align}

The known quantity $\la_i \: \und{F}_{h \restrictto{\Gamma}}^{(i)}$ involves the projection $\Pi_{\cont_h}^{\Gamma}$ of the FE stress field $\cont_h$ over edge $\Gamma \in \Jc_i^{\Ic}$:
\begin{itemize}
\item
for an internal edge $\Gamma$ connected to elements $E$ and $E'$,
\begin{align*}
\displaystyle \Pi_{\cont_h}^{\Gamma} = \frac{1}{2} ( \frac{1}{\eta_{E}} \cont_{h \restrictto{E}} \: \und{n}_E + \frac{1}{\eta_{E'}} \cont_{h \restrictto{E'}} \: \und{n}_{E'} );
\end{align*}
\item
for an edge $\Gamma \subset \partial_1 \Om$ connected to element $E$,
\begin{align*}
\displaystyle \Pi_{\cont_h}^{\Gamma} = \frac{1}{\eta_{E}} \cont_{h \restrictto{E}} \: \und{n}_E;
\end{align*}
\item
for an edge $\Gamma \subset \partial_2 \Om$ connected to element $E$,
\begin{align*}
\displaystyle\Pi_{\cont_h}^{\Gamma} = \frac{1}{\eta_{E}} \und{F}_d.
\end{align*}
\end{itemize}
Then, $\la_i \: \und{F}_{h \restrictto{\Gamma}}^{(i)}$ is defined for each edge $\Gamma \in \Jc_i^{\Ic}$ as:
\begin{itemize}
\item
for $p \geqslant 2$, $\la_i \: \und{F}_{h \restrictto{\Gamma}}^{(i)}$ is the projection of $\la_i \: \Pi_{\cont_h}^{\Gamma}$ over the space of polynomial functions of degree $p$;\\
\item
for $p = 1$, due to implicit conditions (\ref{eq3:wh}), $\la_i \: \und{F}_{h \restrictto{\Gamma}}^{(i)}$ is such that:
\begin{align}\label{eq1:lambda_i_fh_i}
\int_{\Gamma} \la_i \: \und{F}_{h \restrictto{\Gamma}}^{(i)} \: \la_j \dS = \und{0} \quad \forall \: j \in \Ic_{\Gamma}^{\Jc} \: \textrm{and} \: j \neq i; \quad \sum_{j \in \Ic_{\Gamma}^{\Jc}} \la_j \: \und{F}_{h \restrictto{\Gamma}}^{(j)} = \Pi_{\cont_h}^{\Gamma}.
\end{align}
\end{itemize}

It is important to point out that, with respect to the expressions of $\la_i \: \und{F}_h^{(i)}$, those quantities are sought along each edge $\Gamma \in \Jc$ in the same space as $\la_i \: \hat{\und{F}}_h^{(i)}$, namely the space of piecewise continuous polynomial functions of degree $p$.\\

Let us recall that several cost functions will be defined and investigated in \Sect{5}.\\

Computational aspects of the various methods and a detailed analysis of their implementation are addressed in the next section. This section particularly concerns the computation of the local tractions in the EET and EESPT methods.

\section{PRACTICAL IMPLEMENTATION OF THE DIFFERENT METHODS}\label{4}

\subsection{Practical implementation of the EET}\label{4.1}

Local problem $\Pc_i^{\Nc}$ (\ref{eq1:systproj}) associated with a node $i \in \Nc$ takes the matrix form:
\begin{align}\label{eq1:systprojdiscr}
\Bbb^{(i)} \: \hat{\und{b}}^{(i)} = \und{Q}^{(i)},
\end{align}
where $\hat{\und{b}}^{(i)}$ is the unknown generalized vector which contains the combination of vectors $\hat{\und{b}}_{\restrictto{\Gamma}}^{(i)}$ for every edge $\Gamma \in \Jc_i^{\Nc}$, and $\und{Q}^{(i)}$ is the generalized vector which is the combination of vectors $\und{\Qc}_E^{(i)}$ for every element $E \in \Ec_i^{\Nc}$.

Let us recall that the solvability of (\ref{eq1:systprojdiscr}), for both internal node $i \in \Om$ and node $i \in \partial_2 \Om$, requires the verification of a compatibility condition resulting from the FE equilibrium. For such systems, the vectorial equation coming from one, and only one, element $E \in \Ec_i^{\Nc}$ is removed to fix the kernel of matrix $\Bbb^{(i)}$. Thus, system (\ref{eq1:systprojdiscr}) is changed to:
\begin{align}\label{eq1:systprojdiscrblock}
\tilde{\Bbb}^{(i)} \: \hat{\und{b}}^{(i)} = \tilde{\und{Q}}^{(i)}.
\end{align}
Let $n_{ind}^{(i)}$ be the number of independent equations coming from the prolongation condition, \ie the number of rows of matrix $\tilde{\Bbb}^{(i)}$ involved in (\ref{eq1:systprojdiscrblock}). Let $n_{unk}^{(i)}$ be the number of unknowns of system (\ref{eq1:systprojdiscrblock}), \ie the size of vector $\hat{\und{b}}^{(i)}$.

Furthermore, for all nodes $i \in \partial_2 \Om$, conditions $\displaystyle\hat{\und{b}}_{\restrictto{\Gamma}}^{(i)} = \int_{\Gamma} \und{F}_d \: \varphi_i \dS$ over each edge $\Gamma \subset \partial_2 \Om$ are enforced by Lagrange multipliers. Those additional conditions read:
\begin{align}\label{eq1:eqimp}
\Cbb^{(i)} \: \hat{\und{b}}^{(i)} = \und{q}^{(i)} \quad \textrm{for} \ i \in \partial_2 \Om.
\end{align}
Let $n_{enf}^{(i)}$ be the number of enforced equations of system (\ref{eq1:eqimp}), \ie the number of rows of matrix $\Cbb^{(i)}$.

Besides, one performs a minimization step if and only if the following minimization condition holds: $n_{unk}^{(i)} > n_{enf}^{(i)} + n_{ind}^{(i)}$.
In this case, cost function (\ref{eq1:fonctioncoutEET}) can be rewritten mathematically as:
\begin{align}\label{eq1:fonctioncoutEETdiscr}
J(\hat{\und{b}}^{(i)}) = \frac{1}{2} ( \hat{\und{b}}^{(i)} - \und{b}^{(i)} )^{T} \: \Mbb^{(i)} \: ( \hat{\und{b}}^{(i)} - \und{b}^{(i)} ),
\end{align}
where the generalized vector $\und{b}^{(i)}$ is the combination of known vectors $\und{b}_{\restrictto{\Gamma}}^{(i)}$ for every edge $\Gamma \in \Jc_i^{\Nc}$ and $\Mbb^{(i)}$ is a diagonal matrix.\\

Finally, the problem to be solved for each node $i \in \Nc$ depends on the minimization condition:
\begin{itemize}
\item
if $n_{unk}^{(i)} > n_{enf}^{(i)} + n_{ind}^{(i)}$, it consists of minimizing cost function (\ref{eq1:fonctioncoutEETdiscr}) under both constraints (\ref{eq1:eqimp}) enforced over $\partial_2 \Om$ and constraints (\ref{eq1:systprojdiscrblock}) coming from prolongation condition (\ref{eq1:prolongfort}). Thus, introducing the Lagrangian:
\begin{equation}\label{eq1:lagEET}
\begin{aligned}
L(\hat{\und{b}}^{(i)},\und{\Lambda}_{\Cbb}^{(i)},\und{\Lambda}_{\Bbb}^{(i)}) = & \frac{1}{2} ( \hat{\und{b}}^{(i)} - \und{b}^{(i)} )^{T} \: \Mbb^{(i)} \: ( \hat{\und{b}}^{(i)} - \und{b}^{(i)} ) + ( \Cbb^{(i)} \: \hat{\und{b}}^{(i)} - \und{q}^{(i)} )^T \: \und{\Lambda}_{\Cbb}^{(i)} \\
& + ( \tilde{\Bbb}^{(i)} \: \hat{\und{b}}^{(i)} - \tilde{\und{Q}}^{(i)} )^T \: \und{\Lambda}_{\Bbb}^{(i)},
\end{aligned}
\end{equation}
the system to be solved takes the matrix form:
\begin{equation}\label{PbminEET}
\begin{pmatrix}
\Mbb^{(i)} & \Cbb^{(i)^T} & \tilde{\Bbb}^{(i)^T} \\
\Cbb^{(i)} & 0 & 0 \\
\tilde{\Bbb}^{(i)} & 0 & 0 \\
\end{pmatrix}
\begin{bmatrix}
\hat{\und{b}}^{(i)} \\
\und{\Lambda}_{\Cbb}^{(i)} \\
\und{\Lambda}_{\Bbb}^{(i)} \\
\end{bmatrix}
=
\begin{bmatrix}
\Mbb^{(i)} \: \und{b}^{(i)} \\
\und{q}^{(i)} \\
\tilde{\und{Q}}^{(i)} \\
\end{bmatrix},
\end{equation}
where $\und{\Lambda}_{\Cbb}^{(i)}$ and $\und{\Lambda}_{\Bbb}^{(i)}$ correspond to the vectors of Lagrange multipliers associated with constraints (\ref{eq1:eqimp}) and (\ref{eq1:systprojdiscrblock}), respectively;\\
\item
else, it consists of solving explicitly the following system:
\begin{equation}\label{Pbexplicit}
\begin{pmatrix}
\Cbb^{(i)} \\
\tilde{\Bbb}^{(i)} \\
\end{pmatrix}
\begin{bmatrix}
\hat{\und{b}}^{(i)} \\
\end{bmatrix}
=
\begin{bmatrix}
\und{q}^{(i)} \\
\tilde{\und{Q}}^{(i)} \\
\end{bmatrix}.
\end{equation}
It is worthy noticing that, if $n_{unk}^{(i)} < n_{enf}^{(i)} + n_{ind}^{(i)}$, one gives priority to enforced constraints (\ref{eq1:eqimp}) involving traction force density $\und{F}_d$. Thus, in this case, matrix $\tilde{\Bbb}^{(i)}$ and vector $\tilde{\und{Q}}^{(i)}$ are truncated so that the number of constraints coming from (\ref{eq1:systprojdiscrblock}) (\ie resulting from the prolongation condition) is reduced to $n_{unk}^{(i)} - n_{enf}^{(i)}$.
\end{itemize}

Over each edge $\Gamma \in \Jc$, traction $\hat{\und{F}}_{h \restrictto{\Gamma}}$ reads:
\begin{align}\label{eq1:discretisation_fh_hat}
\hat{\und{F}}_{h \restrictto{\Gamma}} = [\varphi_{\restrictto{\Gamma}}]^{T} \: \hat{\und{F}}^{\Gamma},
\end{align}
where $\hat{\und{F}}^{\Gamma}$ is the unknown vector of components of $\hat{\und{F}}_{h \restrictto{\Gamma}}$ over $\Gamma$, and $[\varphi_{\restrictto{\Gamma}}]^{T}$ is the vector of FE shape functions of degree $p$ over $\Gamma$.

Finally, once projections $\hat{\und{b}}^{(i)}$ have been calculated, tractions $\hat{\und{F}}_{h \restrictto{\Gamma}}$ are recovered by solving a linear local system over each edge $\Gamma \in \Jc$ in the form:
\begin{align}\label{eq1:recover_fh_hat}
\Kbb^{\Gamma} \: \hat{\und{F}}^{\Gamma} = \hat{\und{b}}^{\Gamma},
\end{align}
where the generalized vector $\hat{\und{b}}^{\Gamma}$ is the combination of vectors $\hat{\und{b}}_{\restrictto{\Gamma}}^{(i)}$ for each node $i \in \Nc_{\Gamma}^{\Jc}$, and $\Kbb^{\Gamma}$ is similar to a mass matrix associated with edge $\Gamma \in \Jc$.

\subsection{Practical implementation of the SPET}\label{4.2}

In order to simply remove the first (linear) part $\Pi \uu^{\ast}$ of $\uu^{\ast}$ involved in the r.h.s. term of (\ref{eq2:pbpatch}), quantities $\und{e}_i$ and $\uu^{\ast}$ are assumed to belong to the FE interpolation space defined from hierarchical shape functions of degree $p+3$ over each element $E \in \Ec_i^{\Ic}$, instead of Lagrange shape functions. Furthermore, solutions $\und{e}_i$ of (\ref{eq2:pbpatch}) have to verify the following kinematic conditions $\und{e}_{i \restrictto{\partial_1 \Om \cap \Om_i}} = \underline{0}$ (\ie $\und{e}_i \in \Ucb_{0 \restrictto{\Om_i}}$); those kinematic constraints are enforced by substitution in (\ref{eq2:pbpatch}). Therefore, the procedure used for the calculation of solutions $\und{e}_i$ of (\ref{eq2:pbpatch}) is fairly simple to implement as the use of FE hierarchical shape functions enables to retain only the high-order part of $\uu^{\ast}$, which apparently seems to be the main difficulty associated with this method.

\subsection{Practical implementation of the EESPT}\label{4.3}

Searched quantity $\la_i \: \hat{\und{F}}_{h}^{(i)}$ and known quantity $\la_i \: \und{F}_{h}^{(i)}$ are discretized over edge $\Gamma \in \Jc_i^{\Ic}$ in the form:
\begin{align}\label{eq1:discretisation_fh}
\la_i \: \hat{\und{F}}_{h \restrictto{\Gamma}}^{(i)} = [\varphi_{\restrictto{\Gamma}}]^{T} \: \hat{\und{f}}_{h, \Gamma}^{(i)}; \qquad \la_i \: \und{F}_{h \restrictto{\Gamma}}^{(i)} = [\varphi_{\restrictto{\Gamma}}]^{T} \: \und{f}_{h, \Gamma}^{(i)},
\end{align}
where $[\varphi_{\restrictto{\Gamma}}]$ is the matrix of FE shape functions of degree $p$ over $\Gamma$ and $\hat{\und{f}}_{h, \Gamma}^{(i)}$ (respectively $\und{f}_{h, \Gamma}^{(i)}$) is the vector of components of $\la_i \: \hat{\und{F}}_{h \restrictto{\Gamma}}^{(i)}$ (respectively $\la_i \: \und{F}_{h \restrictto{\Gamma}}^{(i)}$).

Besides, quantity $\und{v}_h^{\ast}$ is discretized over element $E \in \Ec_i^{\Ic}$ which reads:
\begin{align}\label{eq1:discretisation_vh}
\und{v}_{h \restrictto{E}}^{\ast} = [\la_{\restrictto{E}}]^{T} \: \und{X}_{h, E}^{(i)},
\end{align}
where $[\la_{\restrictto{E}}]$ is the matrix of linear FE shape functions over $E$ and $\und{X}_{h, E}^{(i)}$ is the vector of components of $\und{v}_{h \restrictto{E}}^{\ast}$.\\

Local problems $\Pc_i^{\Ic}$ (\ref{eq2:pbdensitespatch}) (or (\ref{eq3:pbdensitespatch})) then take the matrix form:
\begin{align}\label{eq1:pblocaldiscr}
\hat{\und{f}}_h^{(i)^{T}} \: \Abb^{(i)} \: \und{X}_h^{(i)} = \und{R}^{(i)^{T}} \: \und{X}_h^{(i)} \quad \forall \: \und{X}_h^{(i)}
\end{align}
where $\hat{\und{f}}_h^{(i)}$ and $\und{X}_h^{(i)}$ are the generalized vectors which contain the combination of vectors $\hat{\und{f}}_{h, \Gamma}^{(i)}$ for every edge $\Gamma \in \Jc_i^{\Ic}$ and vectors $\und{X}_{h, E}^{(i)}$ for every element $E \in \Ec_i^{\Ic}$, respectively. $\Abb^{(i)}$ can be seen as the combination of mass matrices associated with every edge $\Gamma \in \Jc_i^{\Ic}$.\\

Let us recall that conditions $\la_i \: \hat{\und{F}}_{h \restrictto{\Gamma}}^{(i)} = \la_i \: \und{F}_{h \restrictto{\Gamma}}^{(i)}$ along edges $\Gamma \in \Jc_i^{\Ic} \cap \partial_2 \Om$ must be verified for vertex $i \in \partial_2 \Om$. Those equalities read $\hat{\und{f}}_{h, \Gamma}^{(i)} = \und{f}_{h, \Gamma}^{(i)}$ over $\Gamma \in \Jc_i^{\Ic} \cap \partial_2 \Om$ and can be enforced using penalty terms in cost function (\ref{eq1:fonctioncoutEESPT}). Thus, cost function (\ref{eq1:fonctioncoutEESPT}) takes the matrix form:
\begin{align}\label{eq1:fonctioncoutEESPTdiscr}
\frac{1}{2} ( \hat{\und{f}}_{h}^{(i)} - \und{f}_{h}^{(i)} )^{T} \: \Pbb^{(i)} \: ( \hat{\und{f}}_{h}^{(i)} - \und{f}_{h}^{(i)} ),
\end{align}
where the generalized vector $\und{f}_{h}^{(i)}$ is the combination of known vectors $\und{f}_{h, \Gamma}^{(i)}$ for every edge $\Gamma \in \Jc_i^{\Ic}$ and $\Pbb^{(i)}$ is a diagonal penalty matrix.
\begin{itemize}
\item
In the case $p \geqslant 2$, known vector $\und{f}_{h, \Gamma}^{(i)}$ over edge $\Gamma \in \Jc_i^{\Ic}$ is calculated explicitly as it is the projection of $\la_i \: \Pi_{\cont_h}^{\Gamma}$ over the space of polynomial functions of degree $p$, where $\Pi_{\cont_h}^{\Gamma}$ depends on the FE stress field $\cont_h$ and the traction force density $\und{F}_d$.\\
\item
In the case $p = 1$, conditions (\ref{eq1:lambda_i_fh_i}) required for the calculation of quantity $\und{f}_{h, \Gamma}^{(i)}$ over edge $\Gamma \in \Jc_i^{\Ic}$ mathematically read under the form:
\begin{align}\label{eq1:lambda_i_fh_idiscr}
\dbeta_{\Gamma, j}^{(i)} \: \und{f}_{h, \Gamma}^{(i)} = \und{0} \quad \forall \: j \in \Ic_{\Gamma}^{\Jc} \: \textrm{and} \: j \neq i; \quad \sum_{j \in \Ic_{\Gamma}^{\Jc}} \und{f}_{h, \Gamma}^{(j)} = \Pi_{\cont_h}^{\Gamma},
\end{align}
where $\dbeta_{\Gamma, j}^{(i)}$ represents the projection of $[\varphi_{\restrictto{\Gamma}}]$ over edge $\Gamma$ and over linear FE shape function $\la_j$.
Grouping those conditions together, quantities $\und{f}_{h, \Gamma}^{(i)} \quad \forall \: i \in \Ic_{\Gamma}^{\Jc}$ are recovered by solving a linear local system over edge $\Gamma \in \Jc$:
\begin{align}\label{eq1:linearlocalsystem_fh_Gamma}
\Bbb^{\Gamma} \: \und{f}_{h}^{\Gamma} = \und{Q}^{\Gamma},
\end{align}
where $\und{f}_{h}^{\Gamma}$ denotes the generalized vector which is the combination of vectors $\und{f}_{h, \Gamma}^{(i)} \quad \forall \: i \in \Ic_{\Gamma}^{\Jc}$ and $\Bbb^{\Gamma}$ involves matrices $\dbeta_{\Gamma, j}^{(i)} \quad \forall \: j \in \Ic_{\Gamma}^{\Jc} \setminus \left\{ i \right\}, \ \forall \: i \in \Ic_{\Gamma}^{\Jc}$ and identity matrix $\Id$.
\end{itemize}

\begin{figure}
\centering\includegraphics[scale=0.4]{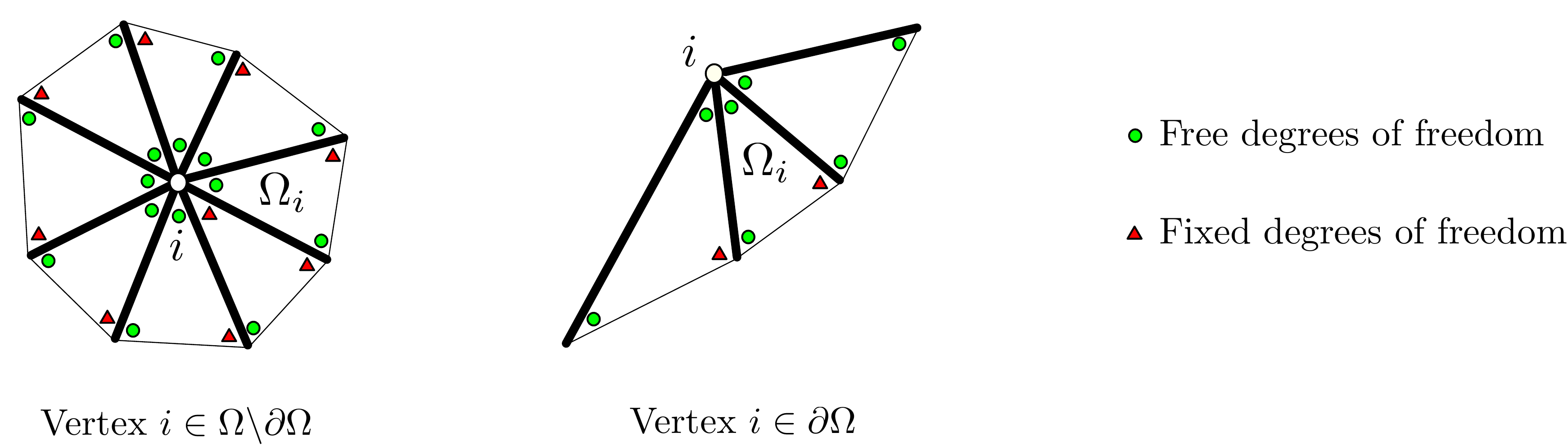}
\caption{Free (green circle) and fixed (red triangle) degrees of freedom over patch $\Om_i$ for an internal vertex $i \in \Om \setminus \dOm$ (left) and a vertex $i \in \dOm$ (right).}\label{fig1:block_kernel_black_modified_with_legend}
\end{figure}

Now, let us recall that kernel $\mathbf{Ker}(\Abb^{(i)})$ is the set of vectors $\und{X}_h^{(i)}$ coming from continuous functions defined over patch $\Om_i$ which vanish on $\Jc_i^{\Ic} \cap \dOm$, \ie $\und{v}_h^{\ast} \in {\bar{\Ucb}_{h,0 \restrictto{\Om_i}}^1}$. In practice, this kernel is automatically reduced to $\und{X}_h^{(i)} = \und{0}$ by fixing some degrees of freedom over patch $\Om_i$. This automatic procedure constitutes a great advantage of the method. \Fig{fig1:block_kernel_black_modified_with_legend} displays free and fixed degrees of freedom for both an internal vertex $i \in \Om \setminus \dOm$ and a vertex $i \in \dOm$.

Once the procedure to fix the kernel of matrix $\Abb^{(i)}$ has been performed, local problem (\ref{eq1:pblocaldiscr}) becomes:
\begin{align}\label{eq2:pblocaldiscr}
\hat{\und{f}}_h^{(i)^{T}} \: \tilde{\Abb}^{(i)} \: \und{\tilde{X}}_h^{(i)} = \tilde{\und{R}}^{{(i)}^{T}} \: \und{\tilde{X}}_h^{(i)} \qquad \forall \: \und{\tilde{X}}_h^{(i)}.
\end{align}

Finally, the problem to be solved over patch $\Om_i$ consists of minimizing cost function (\ref{eq1:fonctioncoutEESPTdiscr}) under constraint (\ref{eq2:pblocaldiscr}) coming from the global prolongation condition (\ref{eq2:prolongfort}). Thus, introducing the Lagrangian:
\begin{align}\label{eq1:lagEESPT}
L(\hat{\und{f}}_h^{(i)},\und{\tilde{X}}_h^{(i)}) = \frac{1}{2} ( \hat{\und{f}}_{h}^{(i)} - \und{f}_{h}^{(i)} )^{T} \: \Pbb^{(i)} \: ( \hat{\und{f}}_{h}^{(i)} - \und{f}_{h}^{(i)} ) + ( \tilde{\Abb}^{(i)^{T}} \: \hat{\und{f}}_h^{(i)} - \tilde{\und{R}}^{(i)} )^T \: \und{\tilde{X}}_h^{(i)},
\end{align}
where $\und{\tilde{X}}_h^{(i)}$ represents the vector of Lagrange multipliers associated with constraints (\ref{eq2:pblocaldiscr}), the system to be solved takes the matrix form:
\begin{equation}\label{PbminEESPT}
\begin{pmatrix}
\Pbb^{(i)} & \tilde{\Abb}^{(i)} \\
\tilde{\Abb}^{(i)^{T}} & 0 \\
\end{pmatrix}
\begin{bmatrix}
\hat{\und{f}}_{h}^{(i)} \\
\und{\tilde{X}}_h^{(i)} \\
\end{bmatrix}
=
\begin{bmatrix}
\Pbb^{(i)} \: \und{f}_{h}^{(i)} \\
\tilde{\und{R}}^{(i)} \\
\end{bmatrix}.
\end{equation}

Finally, referring to section \ref{3.4.2}, one recovers tractions $\hat{\und{F}}_h$ along edges $\Gamma \in \Jc$ from quantities $\hat{\und{f}}_{h, \Gamma}^{(i)}$ for all vertices $i \in \Ic_{\Gamma}^{\Jc}$:
\begin{align}\label{eq1:interpolationFhdiscr}
\hat{\und{F}}_{h \restrictto{\Gamma}} = \sum_{i \in \Ic_{\Gamma}^{\Jc}} [\varphi_{\restrictto{\Gamma}}]^{T} \: \hat{\und{f}}_{h, \Gamma}^{(i)}.
\end{align}

As for the EET, an important and noteworthy point related to the procedure for constructing tractions $\hat{\und{F}}_h = \displaystyle\sum_{i \in \Ic} \la_i \: \hat{\und{F}}_h^{(i)}$ concerns its explicit and non-intrusive nature.

\section{NUMERICAL RESULTS}\label{5}

All the two- and three-dimensional numerical experiments presented in this section are based on finite elements with Lagrange shape functions for both EET and EESPT methods, and hierarchical shape functions for the SPET. All the resolutions of local problems described in \Sect{3} at the element scale for both the EET and EESPT or at the patch scale for the SPET are performed by using a $p + k$ discretization with $k = 3$ ($p$-refinement). The value of the penalty terms involved in the EESPT for enforcing the Neumann conditions on the edges which lie on $\partial_2 \Om$ is set to $10^{5}$. The selected examples are commonly used for industrial applications:
\begin{itemize}
\item
the first case is a weight sensor whose sensor capacity ranges from $1$ kg to $2.5$ kg. This two-dimensional numerical test is a plane-stress linear elastic problem;
\item
the second case is an open hole specimen, which is represented by a three-dimensional plate containing a central hole. This structure is generally employed to study bolted joints in construction and machine design;
\item
the third case is a dome-shaped closure head, forming a part of a nuclear reactor vessel closure assembly;
\item
the fourth case is the hub of the main rotor of the NH90 helicopter developed by the Eurocopter company. It acts as a coupling sleeve between the helicopter frame and the rotor system.

\end{itemize}

For all the structures considered, the material is chosen to be isotropic, homogeneous, linear, and elastic with Young's modulus $E = 1$ and Poisson's ratio $\nu = 0.3$.

The three techniques used to calculate admissible stress fields are analyzed in terms of quality of the error estimators, computational cost and simplicity of practical implementation. The corresponding error estimates are denoted by $\theta = e_{\cre}(\hat{\uu}_h, \hat{\cont}_h)$. The considered error estimates are also compared to the reference error, which usually replaces the unknown exact discretization error for practical purposes. Indeed, the exact error is estimated using a reference error lying in a discretized space much refined with respect to $\Ucb_{h}$; this reference error is thus computed from an ``overkill'' solution. Let us note that the same notation $\lnorm{\und{e}_h}_{u, \Om}$ will be used to denote the energy norm of the reference error and that of the true or exact discretization error. The quality of the error estimators is measured with the usual effectivity index $\eta_{\bullet}$, which is the ratio between any error estimate $\theta_{\bullet}$ and the energy norm of the reference error $\lnorm{\und{e}_h}_{u, \Om}$:
\begin{align}
\eta_{\bullet} = \frac{\theta_{\bullet}}{\lnorm{\und{e}_h}_{u, \Om}}, \notag
\end{align}
the subscript $\bullet$ corresponding to the three estimators studied, \ie \textquotedblleft EET\textquotedblright, \textquotedblleft SPET\textquotedblright \, and \textquotedblleft EESPT\textquotedblright. Thus, the accuracy of the error estimator is given by the value of $\eta$: $\eta > 1$ indicates that $\theta$ is an overestimation of the reference error, whereas $\eta < 1$ demonstrates that the former underestimates the latter.

\subsection{Weight sensor}\label{5.1}

\begin{figure}
\centering\includegraphics[scale = 0.4]{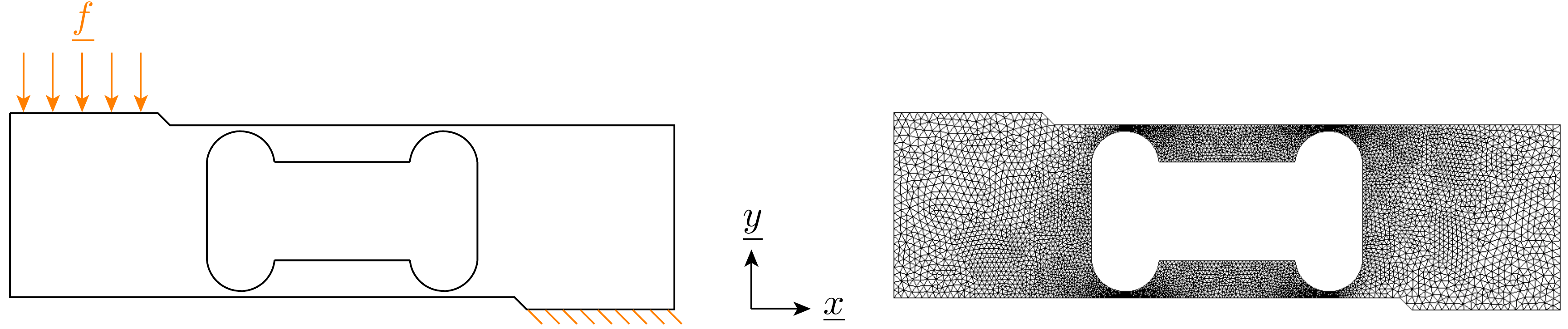}
\caption{Weight sensor model problem (left) and associated finite element mesh (right).}\label{fig1:capteur_effort_2D_geometry_fine_mesh}
\end{figure}

First, let us consider the structure of \Fig{fig1:capteur_effort_2D_geometry_fine_mesh}, which contains two symmetric holes. The structure is clamped along the bottom-right side and subjected to a unit force density $\und{f} = -\und{y}$ along the top-left side. The other remaining boundaries are traction-free. The mesh used for the calculation consists of $11 \, 807$ linear triangular elements and $6 \, 320$ nodes (\ie $12
 \, 640$ degrees of freedom), see \Fig{fig1:capteur_effort_2D_geometry_fine_mesh}. It has been heuristically adapted by refining the elements in the regions of the two holes where stresses are larger (indeed, only the high-stress zones are interesting for design purposes in mechanical design). The reference mesh used to compute the overall reference solution contains $3 \, 326 \, 963$ linear triangular elements and $1 \, 668 \, 711$ nodes (\ie $3 \, 337 \, 422$ d.o.f.). 

\subsubsection{Comparison of the three error estimators}\label{5.1.1}

\begin{figure}
\centering\includegraphics[scale = 0.4]{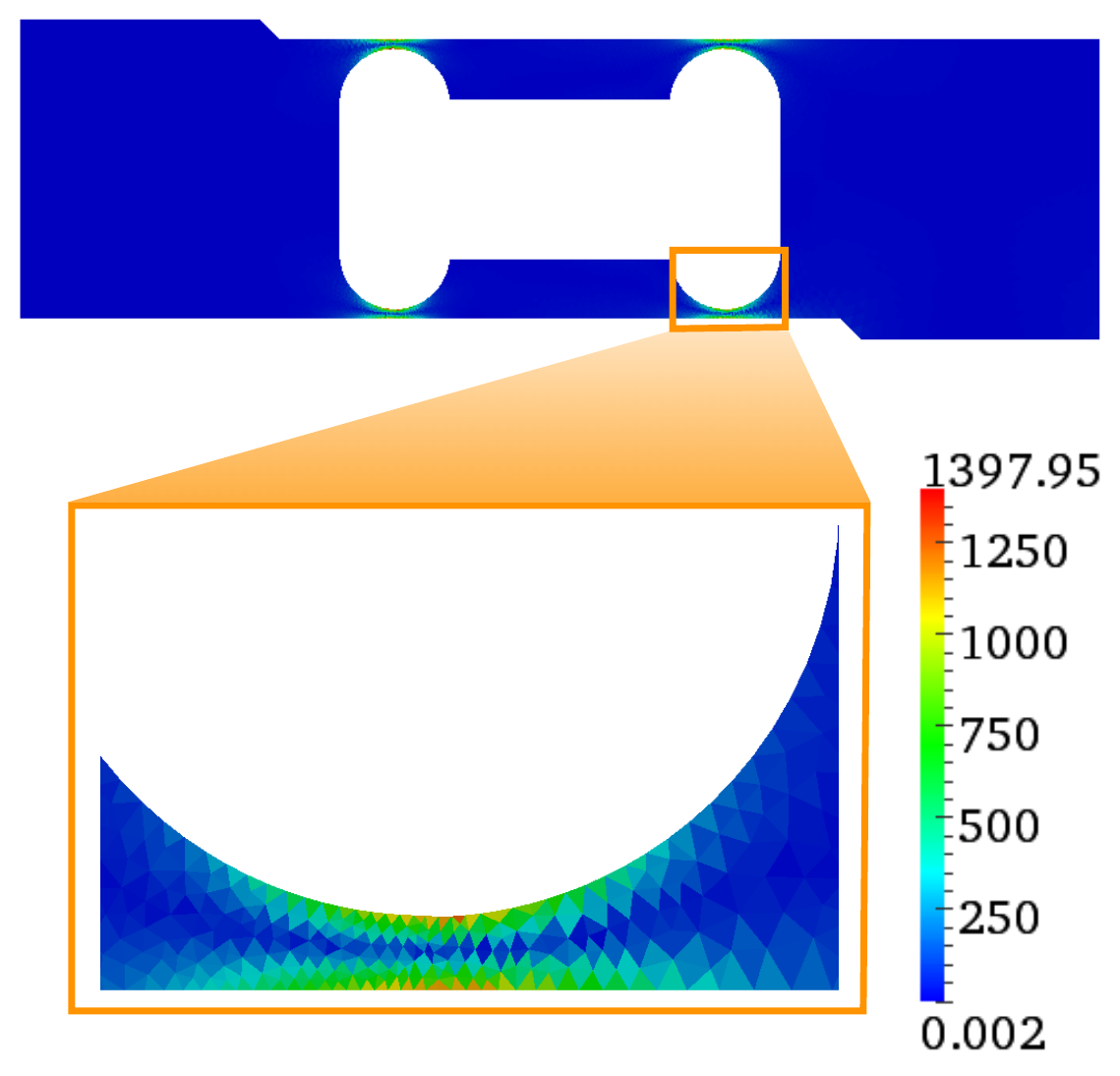}
\caption{Magnitude of the FE stress field and zoom around the highly-loaded region.}\label{fig1:sigma_capteur_effort_2D}
\end{figure}

\begin{figure}
\centering\includegraphics[scale = 0.4]{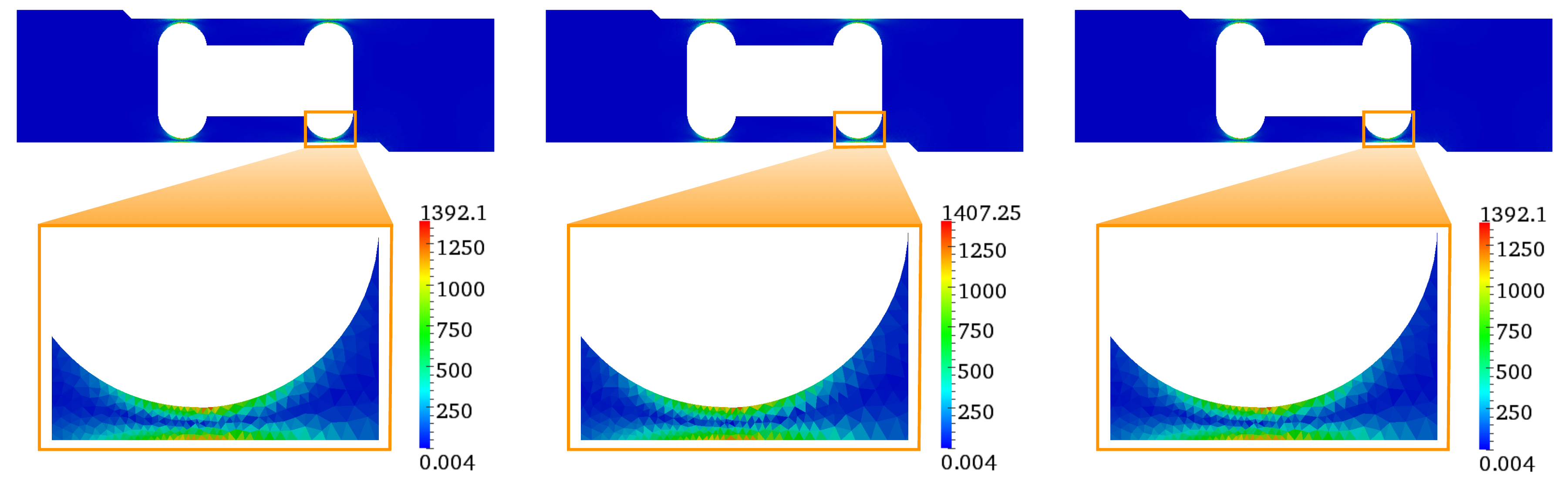}
\caption{Magnitude of the admissible stress field calculated using the EET (left), the SPET (center), and the EESPT (right). Zoom boxes represent the admissible stress fields in the vicinity of the bottom region of the right hole.}\label{fig1:sigma_hat_capteur_effort_2D}
\end{figure}

The local minimization step for both the EET and EESPT is performed using the cost function $J_0$ (\ref{eq1:fonctioncoutEET}), which does not take the element size into account. A detailed analysis of the influence of the choice of the cost function is carried out in the next part. The FE stress field is shown in \Fig{fig1:sigma_capteur_effort_2D} and the admissible stress fields obtained from the three techniques are shown in \Fig{fig1:sigma_hat_capteur_effort_2D}. More precisely, \Fig{fig1:sigma_capteur_effort_2D} and \Fig{fig1:sigma_hat_capteur_effort_2D} display the magnitudes $\displaystyle \sqrt{\Tr\big[\cont_h \: \cont_h \big]}$ and $\displaystyle \sqrt{\Tr\big[\hat{\cont}_h \: \hat{\cont}_h \big]}$ of the FE and admissible stress fields, respectively.

The direct calculation of $\lnorm{\und{e}_h}_{u, \Om}$ is computationally unaffordable. As the structure is only subjected to kinematic conditions imposing zero displacement field on $\partial_1 \Om$, the exact value of the energy norm of the reference error has been calculated using:
\begin{align}
\lnorm{\und{e}_h}_{u, \Om} = \sqrt{\lnorm{\uu}^2_{u, \Om} - \lnorm{\uu_h}^2_{u, \Om}} \simeq \sqrt{\lnorm{\uu_{ref}}^2_{u, \Om} - \lnorm{\uu_h}^2_{u, \Om}} \simeq 347.997.
\end{align}

The required CPU time is about $2$ s, which is very convenient if one seeks to assess the energy norm of the reference error only. Indeed this computation is very low compared to that needed to compute the local contributions to the energy norm of the reference error whose corresponding computational cost reaches about $15$ hours, since the finite element and reference meshes are not being nested and the projection procedure is very time consuming. However, in order to have access to the spatial distribution of the element-by-element contributions to $\lnorm{\und{e}_h}_{u, \Om}$, the calculation of $\lnorm{\und{e}_h}_{u, E} \quad \forall \: E \in \Ec$, with $\lnorm{\bullet}_{u, E} = ( \intE \Tr\big[\K \: \defo(\bullet) \: \defo(\bullet)\big] \dO )^{1/2}$, has been performed.

The performance of the error estimators is summarized in \Tab{table1:comparison_weight_sensor}. First, one verifies that the global effectivity index is higher than $1$ for the three methods. Therefore the three error estimates $\theta_{\bullet}$ are upper bounds with respect to the energy norm of the reference error $\lnorm{\und{e}_h}_{u, \Om}$. The behaviors of the error estimators are consistent. Moreover, the results reveal that the error estimate obtained with the SPET is more accurate than those given by the EET and EESPT: the global effectivity index $\eta$ is closer to $1$ for the SPET. The computational costs required by the EET and the EESPT are quite similar, whereas the CPU time needed for the SPET is four to five times higher. 
\begin{table}
\caption{Comparison of the error estimators given by the EET, the SPET, and the EESPT.}
\centering\small
\begin{tabular}{l c c c c}
\toprule
Methods & Estimate \ $\theta$ & Effectivity index \ $\eta$ & Normalized CPU time \\ 
\midrule
EET & $812.999$ & $2.3362$ & $1.000$ \\ 
SPET & $556.629$ & $1.5995$ & $4.218$ \\ 
EESPT & $812.801$ & $2.3357$ & $1.156$ \\ 
\bottomrule
\end{tabular}
\label{table1:comparison_weight_sensor}
\end{table}

\begin{figure}
\centering\includegraphics[scale = 0.4]{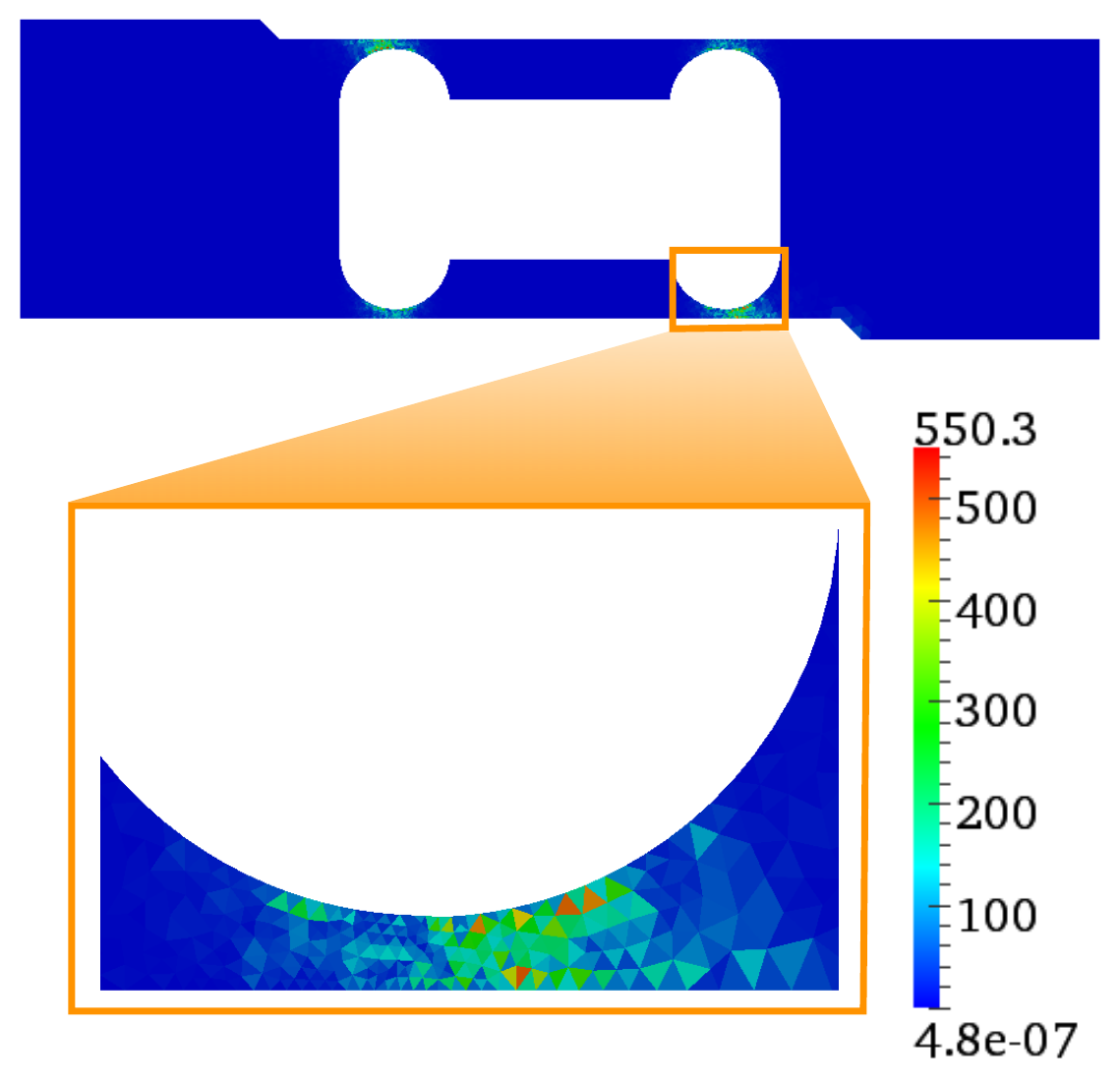}
\caption{Spatial distribution of local contributions to the energy norm of the reference error and zoom around the highly-loaded region.}\label{fig1:distribution_discretization_error_capteur_effort_2D}
\end{figure}

\begin{figure}
\centering\includegraphics[scale = 0.4]{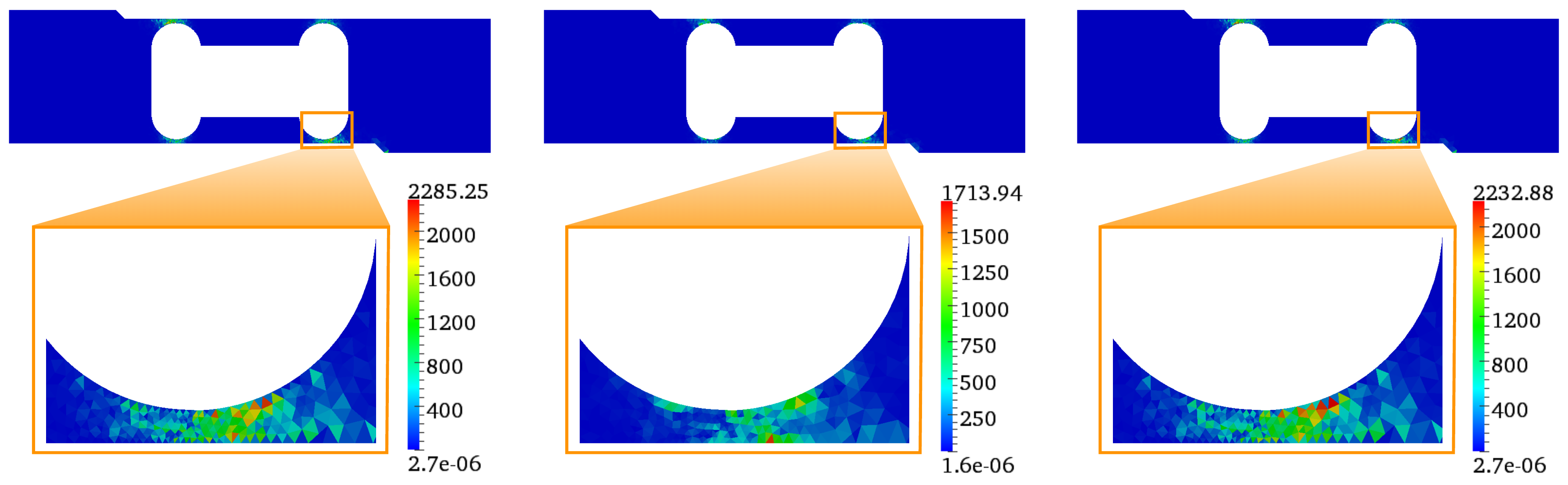}
\caption{Spatial distribution of local contributions to the error estimates computed using the EET (left), the SPET (center), and the EESPT (right). Zoom boxes represent the estimated errors in the vicinity of the bottom region of the right hole.}\label{fig1:distribution_estimators_capteur_effort_2D}
\end{figure}

Spatial distribution of local contributions to the energy norm of the reference error is shown in \Fig{fig1:distribution_discretization_error_capteur_effort_2D}, whereas that of local contributions to the error estimates for the three methods are displayed in \Fig{fig1:distribution_estimators_capteur_effort_2D}. Local contributions in each element of the FE mesh to the error estimates provide local indicators of the local energy norm of the reference error, that could be useful for a remeshing strategy in an adaptive procedure. From a qualitative viewpoint, the estimated error distribution is in good agreement with the reference error distribution. The main contributions of the error are concentrated around the top and bottom regions of the holes, \ie located in the highly-loaded regions.

\begin{figure}
\centering\includegraphics[scale = 0.4]{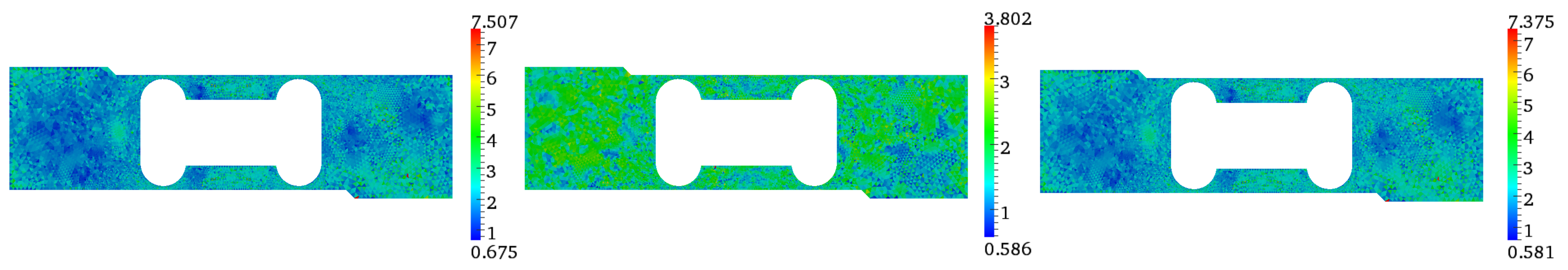}
\caption{Spatial distribution of the local effectivity indices calculated using the EET (left), the SPET (center), and the EESPT (right).}\label{fig1:i_eff_estimators_capteur_effort_2D}
\end{figure}

\newpage
One defines a local effectivity index as the ratio between the element-by-element contribution to the error estimate and the one to the energy norm of the reference error. Spatial distribution of the local effectivity indices is shown in \Fig{fig1:i_eff_estimators_capteur_effort_2D}. Let us recall that contributions to the effectivity index are computed with respect to the reference error. The local effectivity indices range between $0.675$ and $7.507$ for the EET; $0.586$ and $3.802$ for the SPET; $0.581$ and $7.375$ for the EESPT. One can observe that a quite large set of elements has a local effectivity index inferior to $1$, even in low-stress regions.

\subsubsection{Influence of the choice of the cost function involved in the EET and the EESPT}\label{5.1.2}

Quantities ${ \Delta \und{q}^{(i)}}_{\restrictto{\Gamma}}$ involved in the cost function are $( \hat{\und{b}}^{(i)} - \und{b}^{(i)} )_{\restrictto{\Gamma}}$ and $( \la_i \: \hat{\und{F}}_h^{(i)} - \la_i \: \und{F}_h^{(i)} )_{\restrictto{\Gamma}}$ for the EET and the EESPT, respectively. Several cost functions have been considered for the least-squares minimization step over $\Jc_i^{\square}$, where {\scriptsize$\square$} stands for $\Nc$ (resp. $\Ic$) in the case of the EET (resp. EESPT):

\begin{enumerate} 
\item[(i)]
$J_0(\und{q}^{(i)}) = \frac{1}{2} \displaystyle \sum_{\Gamma \in \Jc_i^{\square}} ( \Delta \und{q}^{(i)} \cdot \Delta \und{q}^{(i)} )_{\restrictto{\Gamma}}$; 
\item[(ii)]
$J_1(\und{q}^{(i)}) = \frac{1}{2} \displaystyle \sum_{\Gamma \in \Jc_i^{\square}} \frac{1}{l_{\Gamma}^2} \: ( \Delta \und{q}^{(i)} \cdot \Delta \und{q}^{(i)} )_{\restrictto{\Gamma}}$ which takes the size $l_{\Gamma}$ of each edge $\Gamma \in \Jc_i^{\square}$ into account;
\item[(iii)]
$J_2(\und{q}^{(i)}) = \frac{1}{2} \displaystyle \sum_{\Gamma \in \Jc_i^{\square}} \frac{1}{l_{\Gamma}^2} \: \frac{1 + \nu}{E} \left( \frac{1 - 2 \nu}{1 - \nu} (\Delta \und{q}^{(i)} \cdot \und{n})_{\restrictto{\Gamma}}^2 + 2 \: (\Pi \Delta \und{q}^{(i)} \cdot \Pi \Delta \und{q}^{(i)})_{\restrictto{\Gamma}} \right)$ which corresponds to the density of elastic energy stored in each edge $\Gamma \in \Jc_i^{\square}$. It has been explicitly obtained by assuming that the plane part of the strain field associated with the statically admissible stress field is continuous across the inter-element edges. The interested reader can refer to \cite{Flo02bis} for details.
\end{enumerate}

Their influence on the quality of the two error estimators EET and EESPT has been studied; results are shown in \Tab{table1:influence_cost_function_weight_sensor}. They reveal that the more physical information the cost function contains, the more accurate the yielded error estimator is. Indeed the global effectivity index experiences a $8 \%$ decrease when the cost function $J_2$ is used instead of $J_0$. The bounds of the local effectivity indices calculated using each cost function are shown in \Tab{table1:influence_cost_function_i_eff_weight_sensor} for both EET and EESPT methods.
\begin{table}
\caption{Influence of the choice of the cost function involved in the EET and the EESPT on the quality of the error estimators.}
\centering\small
\begin{tabular}{l p{0.2cm} c c c p{0.2cm} c c c}
\toprule
 & & \multicolumn{3}{c}{Estimate \ $\theta$} & & \multicolumn{3}{c}{Effectivity index \ $\eta$} \\
\ms \cline{3-5} \cline{7-9} \ms
Methods & & $J_0$ & $J_1$ &$J_2$ & & $J_0$ & $J_1$ & $J_2$ \\
\midrule
EET & & $812.999$ & $812.801$ & $749.732$ & & $2.3362$ & $2.3357$ & $2.1544$ \\
EESPT & & $812.801$ & $815.754$ & $754.667$ & & $2.3357$ & $2.3441$ & $2.1685$ \\
\bottomrule
\end{tabular}
\label{table1:influence_cost_function_weight_sensor}
\end{table}
\begin{table}
\caption{Influence of the choice of the cost function involved in the EET and the EESPT on the bounds of the local effectivity indices.}
\centering\small
\begin{tabular}{l p{0.2cm} c c c}
\toprule
 & & \multicolumn{3}{c}{Bounds of the local effectivity indices \ $\eta_E$} \\ 
\ms \cline{3-5} \ms
 Methods & & $J_0$ & $J_1$ &$J_2$ \\
\midrule
EET & & $0.675$ - $7.507$ & $0.581$ - $7.375$ & $0.639$ - $7.305$ \\
EESPT & & $0.581$ - $7.375$ & $0.570$ - $7.209$ & $0.610$ - $7.151$ \\
\bottomrule
\end{tabular}
\label{table1:influence_cost_function_i_eff_weight_sensor}
\end{table}

These results confirm that it seems better to use a physically sound cost function, such as $J_2$, even though the bounds of the local effectivity index are quite similar for all the cost functions we considered.
 
\subsection{Plate with a hole}\label{5.2}

\begin{figure}
\centering\includegraphics[scale = 0.4]{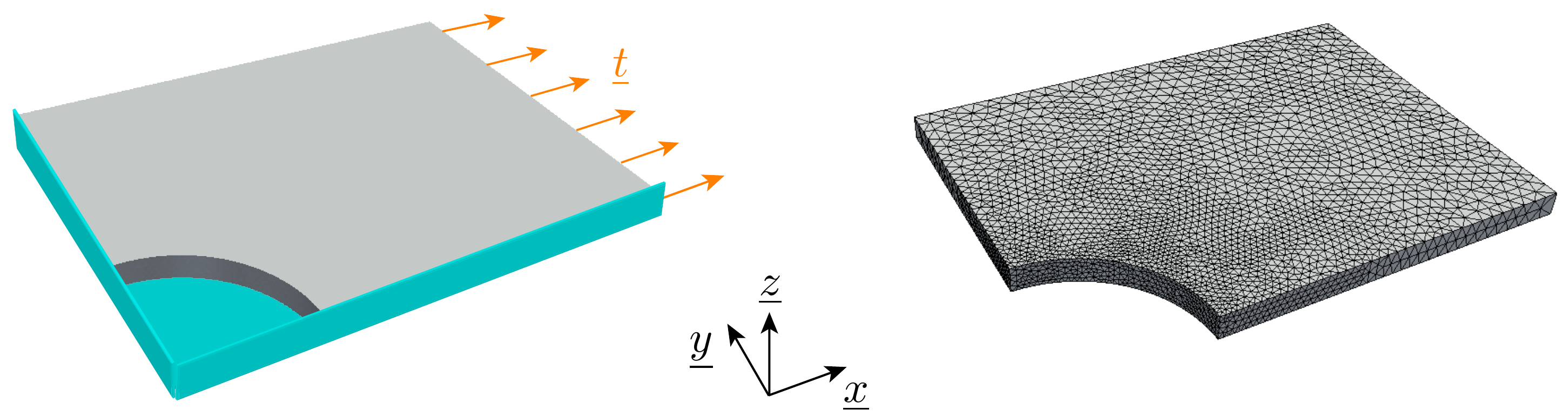}
\caption{Plate with a hole model problem (left) and associated finite element mesh (right).}\label{fig1:plaque_trouee_3D_geometry_fine_mesh}
\end{figure}

Let us consider a plate, represented in \Fig{fig1:plaque_trouee_3D_geometry_fine_mesh}, which contains a central hole. The plate is $20$ mm long, $15$ mm large, $1$ mm high and presents a hole of radius $2.5$ mm. Due to symmetry, only one eighth of the structure is analyzed. Symmetry boundary conditions are applied on the light blue surfaces represented in \Fig{fig1:plaque_trouee_3D_geometry_fine_mesh}. The structure is subjected to a unit traction force density $\und{t} = +\und{x}$ along the right side. The hole and the top side are traction-free boundaries. The mesh, containing $23 \, 493$ linear tetrahedral elements and $6 \, 125$ nodes (\ie $18 \, 375$ d.o.f.), is given in \Fig{fig1:plaque_trouee_3D_geometry_fine_mesh}. The mesh density increases toward the hole which is the highest stress region. The reference mesh is built up by splitting each tetrahedron into $64$ tetrahedra. Consequently, it consists of $1 \, 503\, 552$ linear tetrahedral elements and $284 \, 753$ nodes (\ie $854 \, 259$ d.o.f.).

\subsubsection{Comparison of the three error estimators}\label{5.2.1}

\begin{figure}
\centering\includegraphics[scale = 0.4]{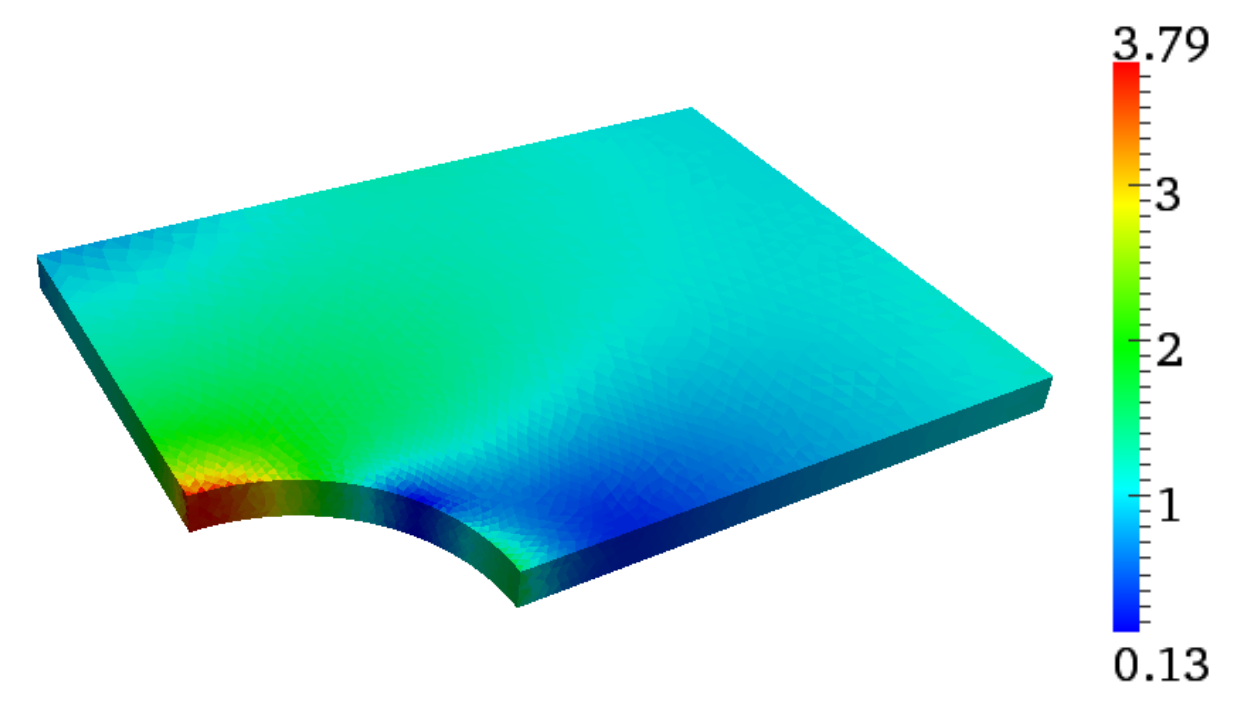}
\caption{Magnitude of the FE stress field.}\label{fig1:sigma_plaque_trouee_3D}
\end{figure}

\begin{figure}
\centering\includegraphics[scale = 0.4]{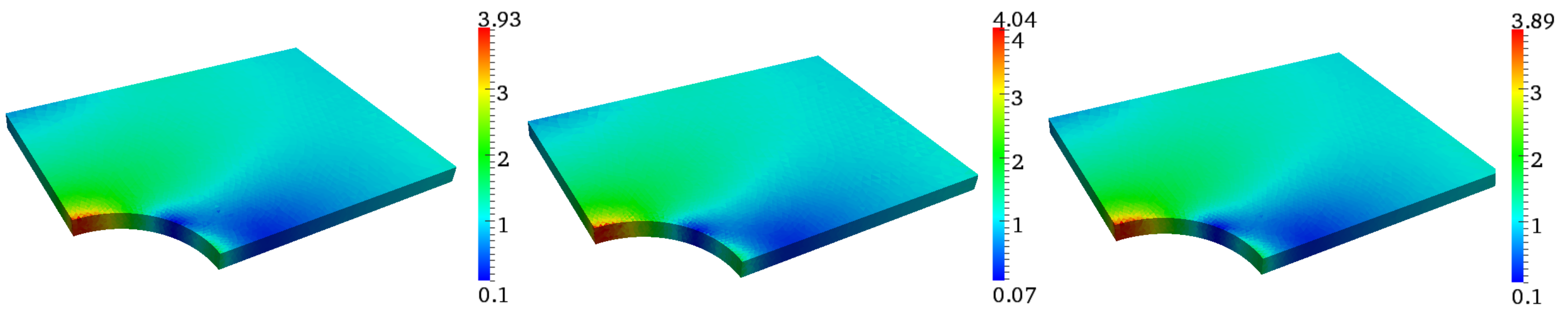}
\caption{Magnitude of the admissible stress field calculated using the EET (left), the SPET (center), and the EESPT (right).}\label{fig1:sigma_hat_plaque_trouee_3D}
\end{figure}

The cost function $J_0$ has been used for the local minimization step. The FE stress field is shown in \Fig{fig1:sigma_plaque_trouee_3D} and the admissible stress fields obtained from the three techniques are displayed in \Fig{fig1:sigma_hat_plaque_trouee_3D}. The highly-loaded region is located toward the vicinity of the hole.

The exact value of the energy norm of the reference error has been directly calculated from the reference solution:
\begin{align}
\lnorm{\und{e}_h}_{u, \Om} = \sqrt{\lnorm{\uu}^2_{u, \Om} - \lnorm{\uu_h}^2_{u, \Om}} \simeq 0.1688.
\end{align}

The CPU time required for this last calculation is about $1$ s, which is very low compared to that needed to compute the element-by-element contributions to $\lnorm{\und{e}_h}_{u, \Om}$ whose computational cost exceeds $22$ hours.

For each method, the estimate, the corresponding global effectivity index and the normalized CPU time (with respect to that required for EET) have been calculated and compared; results are given in \Tab{table1:comparison_plate_with_a_hole}. The same conclusions as for the previous two-dimensional case can be drawn, that is to say that the SPET is more accurate and costly than the two other estimators. Indeed, the computational cost required to compute the SPET is about five and a half times higher than that needed to compute the EET and the EESPT.
\begin{table}
\caption{Comparison of the error estimators given by the EET, the SPET, and the EESPT.}
\centering\small
\begin{tabular}{l c c c c}
\toprule
Methods & Estimate \ $\theta$ & Effectivity index \ $\eta$ & Normalized CPU time \\ 
\midrule
EET & $0.9227$ & $5.4671$ & $1.000$ \\ 
SPET & $0.5928$ & $3.5127$ & $5.368$ \\ 
EESPT & $0.9119$ & $5.4031$ & $1.009$ \\ 
\bottomrule
\end{tabular}
\label{table1:comparison_plate_with_a_hole}
\end{table}

\begin{figure}
\centering\includegraphics[scale = 0.4]{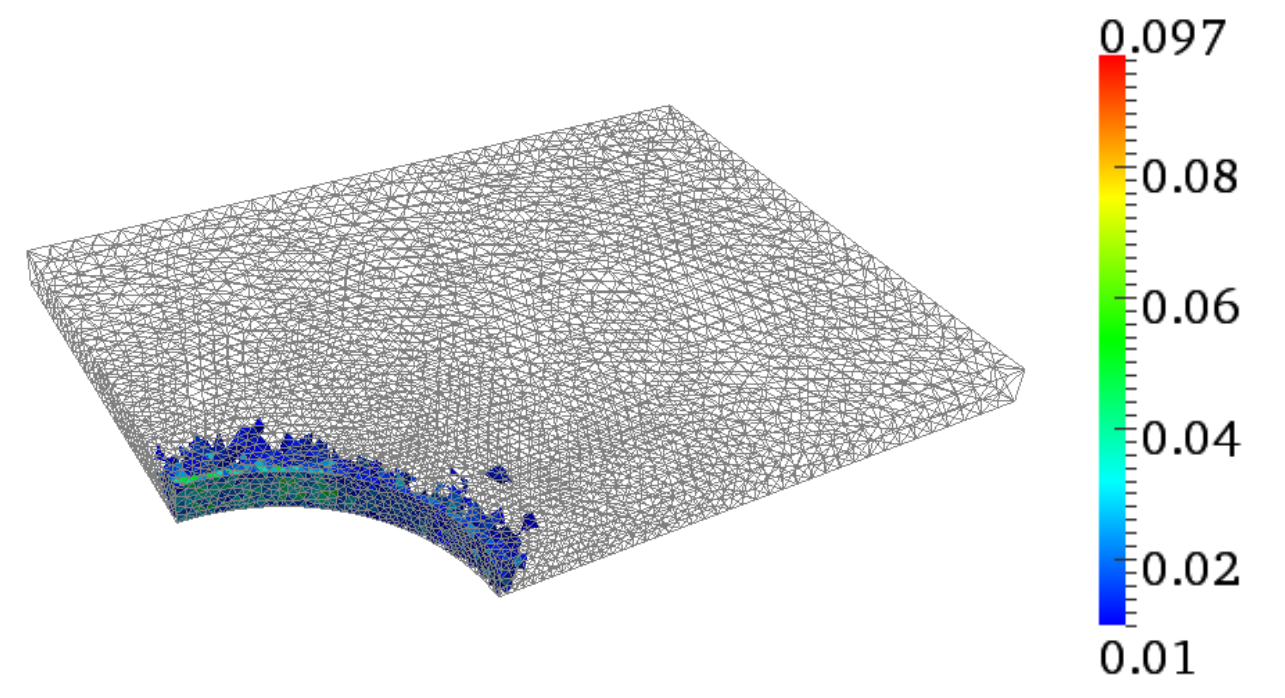}
\caption{Spatial distribution of relevant local contributions to the density of the energy norm of the reference error.}\label{fig1:distribution_discretization_error_plaque_trouee_3D}
\end{figure}

\begin{figure}
\centering\includegraphics[scale = 0.4]{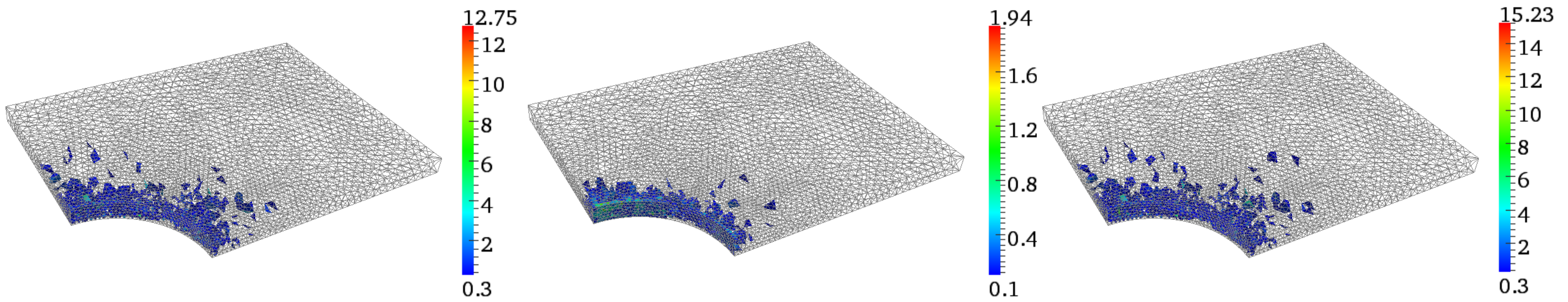}
\caption{Spatial distribution of relevant local contributions to the density of the error estimates calculated using the EET (left), the SPET (center), and the EESPT (right).}\label{fig1:distribution_estimators_plaque_trouee_3D}
\end{figure}

As the density of the FE mesh is not uniform, the reference and estimated spatial distributions of the error are represented as a density that is the ratio between the squared element-by-element contribution to the reference (or estimated) error and the size of the element.
\Figs{fig1:distribution_discretization_error_plaque_trouee_3D} and \ref{fig1:distribution_estimators_plaque_trouee_3D} show respectively the spatial distribution of the local contributions to the density of the energy norm of the reference error and that of the local contributions to the density of the error estimates for the three methods. Only elements in which the contribution to the error estimates is relevant are represented. The main contributions of the density of the error are related to elements located in the neighborhood of the hole. The similarity of maps in \Fig{fig1:distribution_estimators_plaque_trouee_3D} demonstrates the good agreement between the reference and estimated error distributions. Nevertheless, the SPET produces local contributions of the density of the error estimate of better quality level than that given by the EET and the EESPT.


Once again, the local contributions to the effectivity index are computed with respect to the reference error. The local effectivity indices range between $1.253$ and $34.930$ for the EET; $1.086$ and $9.458$ for the SPET; $0.999$ and $34.377$ for the EESPT. One can observe that almost all the elements have a local effectivity index superior to $1$ for the three techniques. Nevertheless, it is worthy noticing that the introduced global error estimators could not directly be used in a goal-oriented analysis, especially the EET and EESPT, because the maximum local contribution to the error estimates drastically overestimates the corresponding contribution to the energy norm of the reference error.

\subsubsection{Influence of the choice of the cost function involved in the EET and the EESPT}\label{5.2.2}

The study of the influence of the various cost functions on the performance of the two error estimators EET and EESPT has been investigated; results are summarized in \Tab{table1:influence_cost_function_plate_with_a_hole}. They tend to show that the global effectivity index experiences a moderate decrease by $3 \%$ for the EET and a slight increase by $1 \%$ for the EESPT between the cost functions $J_0$ and $J_2$. The use of the cost function $J_2$ enables to obtain a slightly better effectivity index compared to the cost function $J_1$ for both estimators. 
\begin{table}
\caption{Influence of the choice of the cost function involved in the EET and the EESPT on the quality of the error estimators.}
\centering\small
\begin{tabular}{l p{0.2cm} c c c p{0.2cm} c c c}
\toprule
 & & \multicolumn{3}{c}{Estimate \ $\theta$} & & \multicolumn{3}{c}{Effectivity index \ $\eta$} \\
\ms \cline{3-5} \cline{7-9} \ms
Methods & & $J_0$ & $J_1$ &$J_2$ & & $J_0$ & $J_1$ & $J_2$ \\
\midrule
EET & & $0.9227$ & $0.9119$ & $0.8906$ & & $5.4671$ & $5.4031$ & $5.2769$ \\
EESPT & & $0.9119$ & $0.9404$ & $0.9136$ & & $5.4031$ & $5.5720$ & $5.4134$ \\
\bottomrule
\end{tabular}
\label{table1:influence_cost_function_plate_with_a_hole}
\end{table}

\begin{table}
\caption{Influence of the choice of the cost function involved in the EET and the EESPT on the bounds of the local effectivity indices.}
\centering\small
\begin{tabular}{l p{0.2cm} c c c}
\toprule
 & & \multicolumn{3}{c}{Bounds of the local effectivity indices \ $\eta_E$} \\ 
\ms \cline{3-5} \ms
 Methods & & $J_0$ & $J_1$ &$J_2$ \\
\midrule
EET & & $1.253$ - $34.930$ & $0.999$ - $34.377$ & $1.076$ - $40.694$ \\
EESPT & & $0.999$ - $34.377$ & $0.885$ - $35.401$ & $0.785$ - $42.717$ \\
\bottomrule
\end{tabular}
\label{table1:influence_cost_function_i_eff_plate_with_a_hole}
\end{table}

The bounds of the local effectivity indices obtained using each cost function are given in \Tab{table1:influence_cost_function_i_eff_plate_with_a_hole} for both the EET and EESPT. The analysis of the data reveals that the two global estimators are not well-suited to direct goal-oriented analysis for all the cost functions considered.

\subsection{Closure head}\label{5.3}

\begin{figure}
\centering\includegraphics[scale = 0.4]{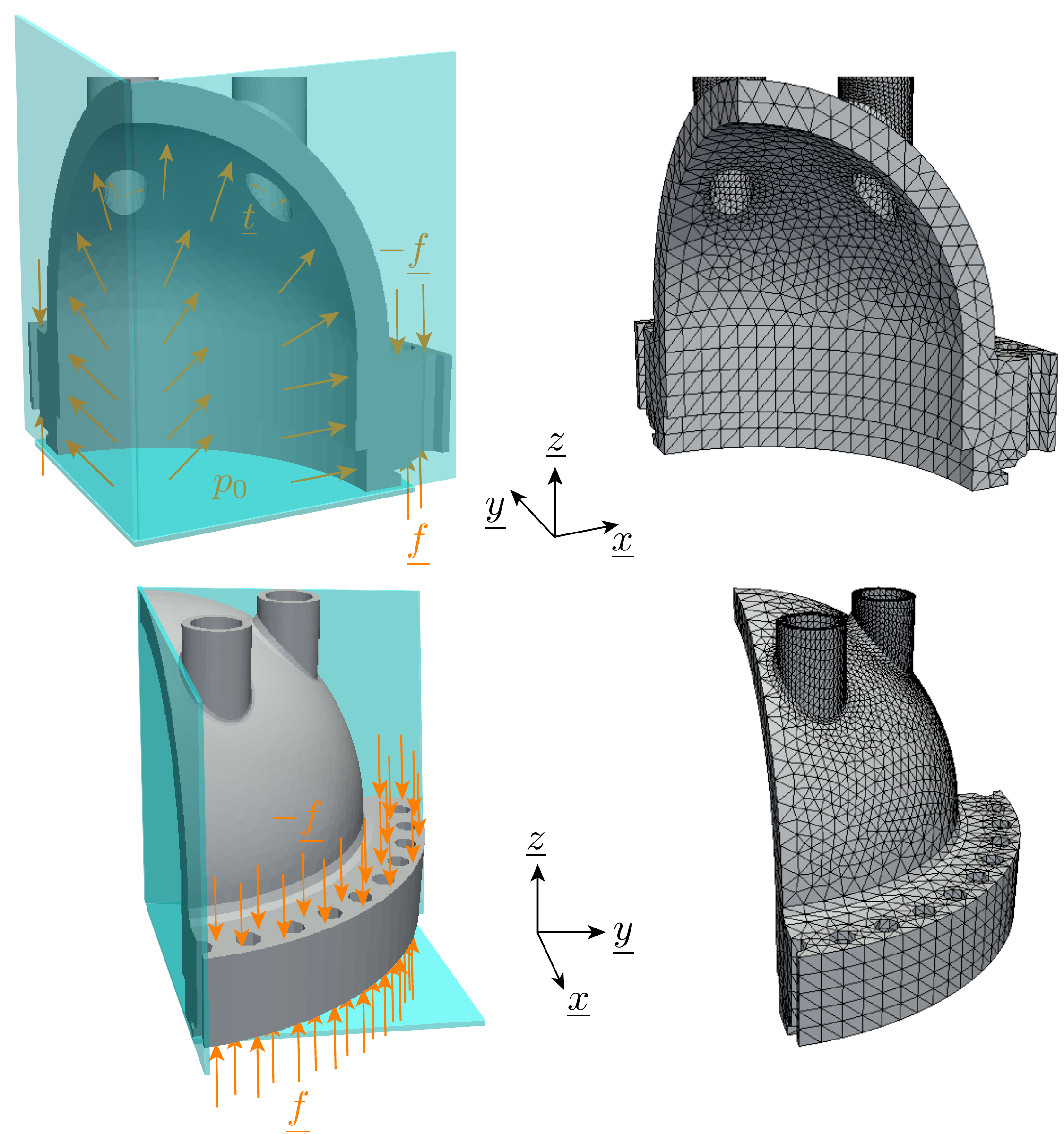}
\caption{Closure head model problem (left) and associated finite element mesh (right). Light blue plans represent symmetry conditions.}\label{fig1:tete_reacteur_3D_geometry_mesh}
\end{figure}

Now, let us consider a dome-shaped closure head which is part of a nuclear reactor vessel closure assembly. This closure head includes height standpipes. In the whole vessel assembly, a control rod drive mechanism (CDM) plug is inserted in each standpipe. The closure head is cyclic-symmetric with respect to the axis of the closure head and only one quarter of the structure is modeled (see \Fig{fig1:tete_reacteur_3D_geometry_mesh}). The structure is $H = 2330$ mm high with an inner radius of $1310$ mm and a thickness of $210$ mm. The bottom of the closure head presents a shoulder that is $70$ mm long, representing its connection to the vessel shell. The standpipes have a height of about $530$ mm above the closure head, an inside diameter of $305$ mm and an outside diameter of $406$ mm. Each standpipe presents a $260$ mm diameter ledge in which a CDM plug is sit in the whole vessel assembly. The center of each standpipe is located at a distance of $964$ mm from the axis of the closure head. Furthermore, the closure head includes a bolting flange containing a set of holes whose diameter is $146$ mm. In the whole closure assembly, the closure head is attached to the vessel shell by $40$ stud bolts passing through this flange, which has a height of $584$ mm.\\

Symmetry boundary conditions are introduced on the bottom end of the closure head and on the perpendicular sides. The inner surface of the closure head is subjected to a constant pressure $p_0 = 1$. A unit traction force density $\und{t} = - \und{n}$, normal to the surface, is applied on the inner surface of the standpipe located under the ledge. Those two loading conditions tend to represent the effect of water on the closure head. Besides, both sides of the flange are loaded with a vertical unit traction force density $\und{f} = \pm \und{z}$, representing the effect of the pre-tension load applied to each bolt on the flange. The geometry and mesh considered, consisting of $38 \, 099$ linear tetrahedral elements and $8 \, 730$ nodes (\ie $26 \, 190$ d.o.f.), are shown in \Fig{fig1:tete_reacteur_3D_geometry_mesh}. The reference mesh is obtained dividing each tetrahedron in $8$ tetrahedra. As a result, it comprises $304 \, 792$ linear tetrahedral elements and $60 \, 381$ nodes (\ie $181 \, 143$ d.o.f.).

\subsubsection{Comparison of the three error estimators}\label{5.3.1}

\begin{figure}
\centering\includegraphics[scale = 0.4]{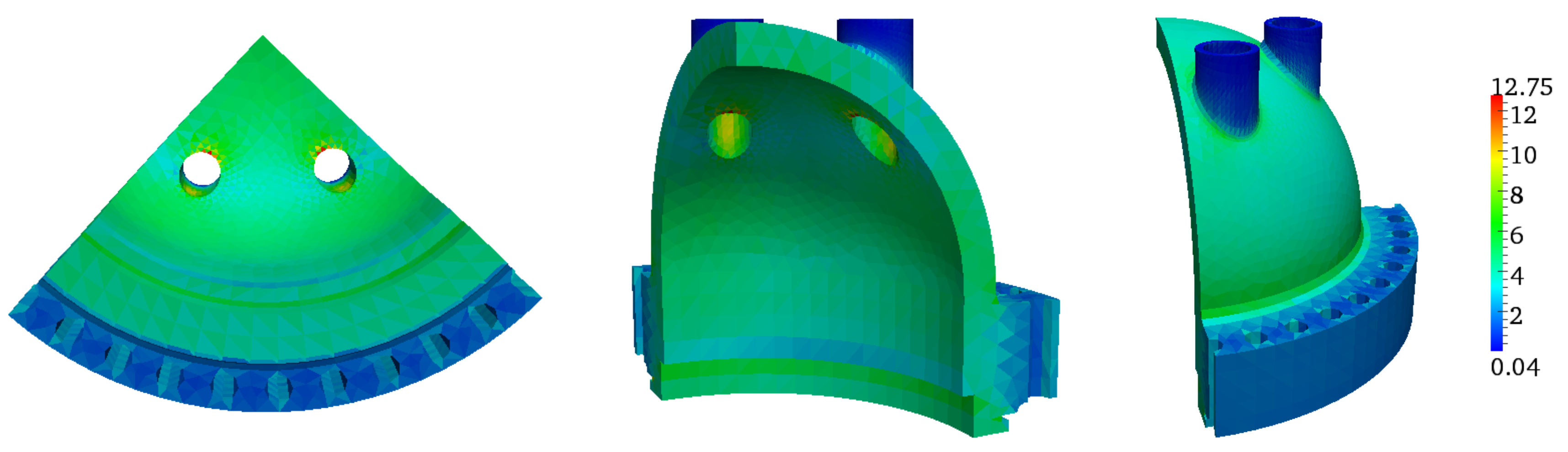}
\caption{Magnitude of the FE stress field.}\label{fig1:sigma_tete_reacteur_3D}
\end{figure}

\begin{figure}
\centering\includegraphics[scale = 0.4]{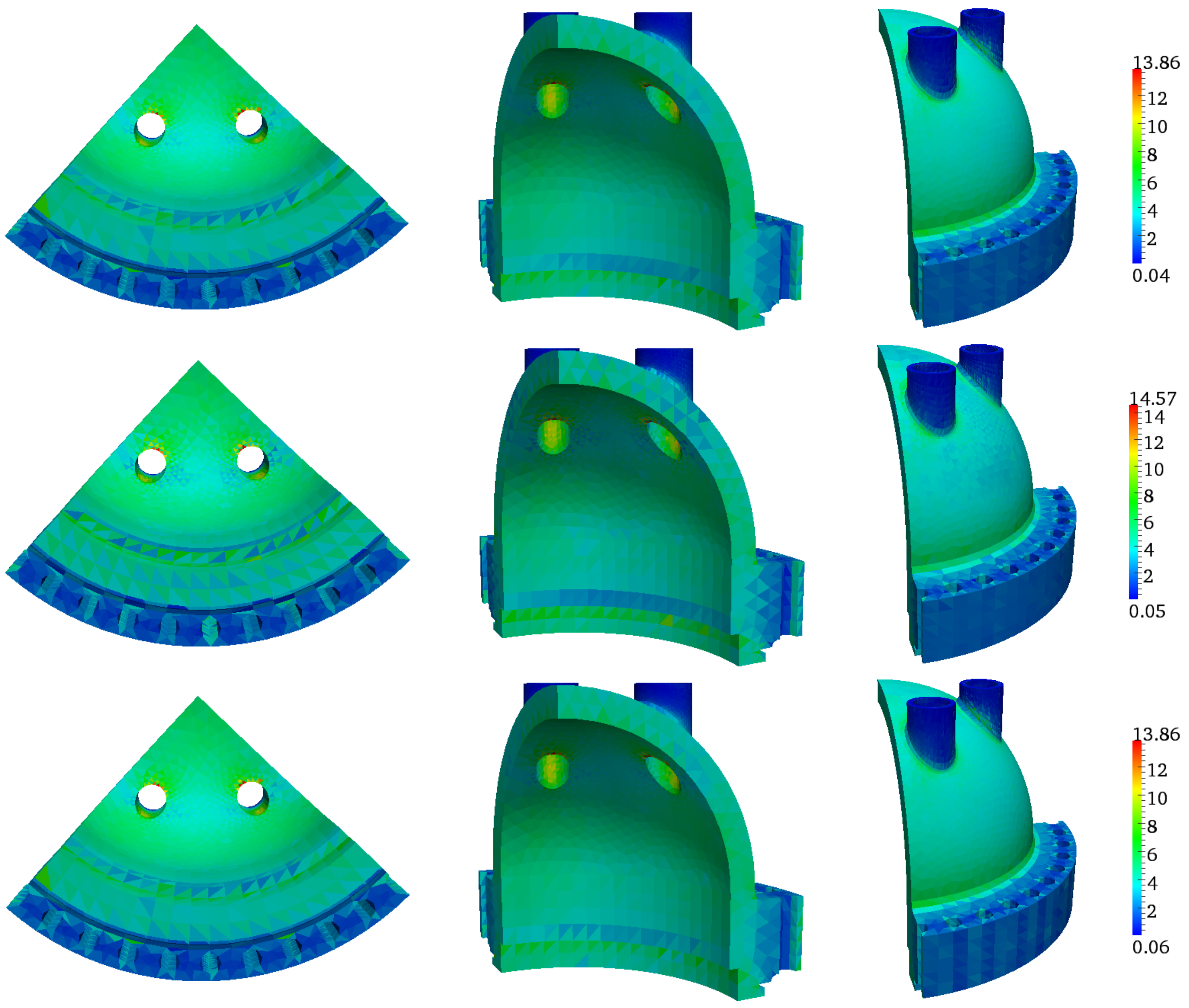}
\caption{Magnitude of the admissible stress field calculated using the EET (top), the SPET (middle), and the EESPT (bottom).}\label{fig1:sigma_hat_tete_reacteur_3D}
\end{figure}

\Figs{fig1:sigma_tete_reacteur_3D} and \ref{fig1:sigma_hat_tete_reacteur_3D} represent the FE stress field and the admissible stress fields, respectively, obtained using the three techniques. A good match is observed for the different stress distributions on the closure head. Results reveal that the highest stress zone is located at the bottom region of the standpipes.

The exact value of the energy norm of the reference error has been directly calculated from the reference solution:
\begin{align}
\lnorm{\und{e}_h}_{u, \Om} = \sqrt{\lnorm{\uu}^2_{u, \Om} - \lnorm{\uu_h}^2_{u, \Om}} \simeq 135.434,
\end{align}
and required a computational cost of about $0.25$ s, whereas the computational cost necessary for calculating the local contributions to $\lnorm{\und{e}_h}_{u, \Om}$ reaches $7$ hours.

The effectiveness and the normalized CPU time corresponding to each estimator are compared again. \Tab{table1:comparison_closure_head} shows the estimate, the effectivity index and the normalized computational cost for the three estimators considered. Those results confirm the relevance and the consistency of the proposed estimators.

\begin{table}
\caption{Comparison of the error estimators given by the EET, the SPET, and the EESPT.}
\centering\small
\begin{tabular}{l c c c c}
\toprule
Methods & Estimate \ $\theta$ & Effectivity index \ $\eta$ & Normalized CPU time \\ 
\midrule
EET & $812.771$ & $6.001$ & $1.000$ \\ 
SPET & $604.311$ & $4.462$ & $5.017$ \\ 
EESPT & $787.983$ & $5.818$ & $1.007$ \\ 
\bottomrule
\end{tabular}
\label{table1:comparison_closure_head}
\end{table}

\begin{figure}
\centering\includegraphics[scale = 0.4]{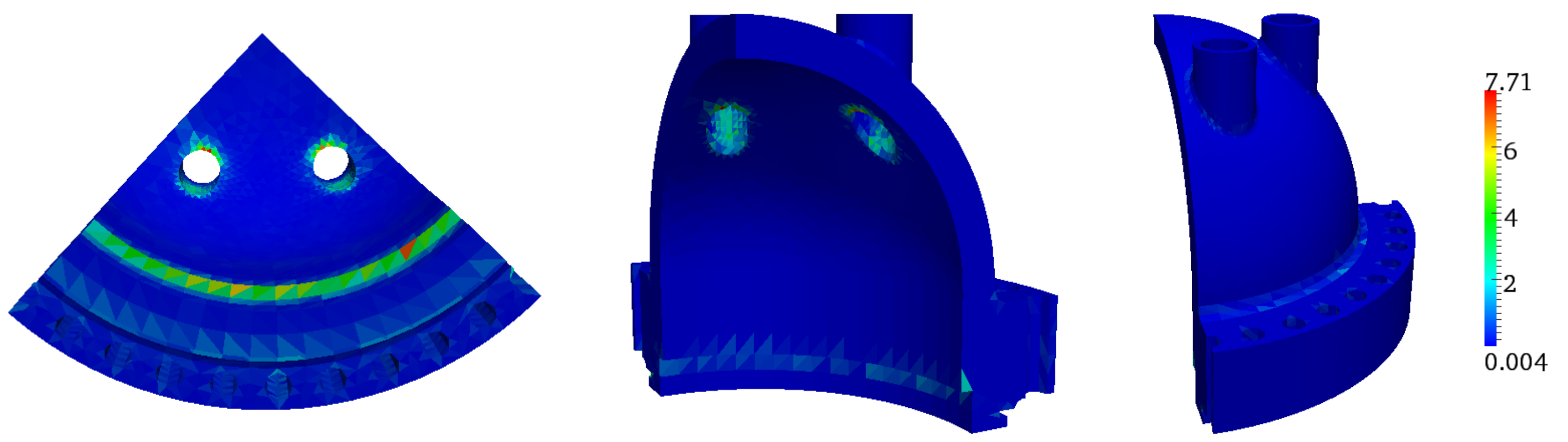}
\caption{Spatial distribution of local contributions to the density of the energy norm of the reference error.}\label{fig1:distribution_discretization_error_tete_reacteur_3D}
\end{figure}

\begin{figure}
\centering\includegraphics[scale = 0.4]{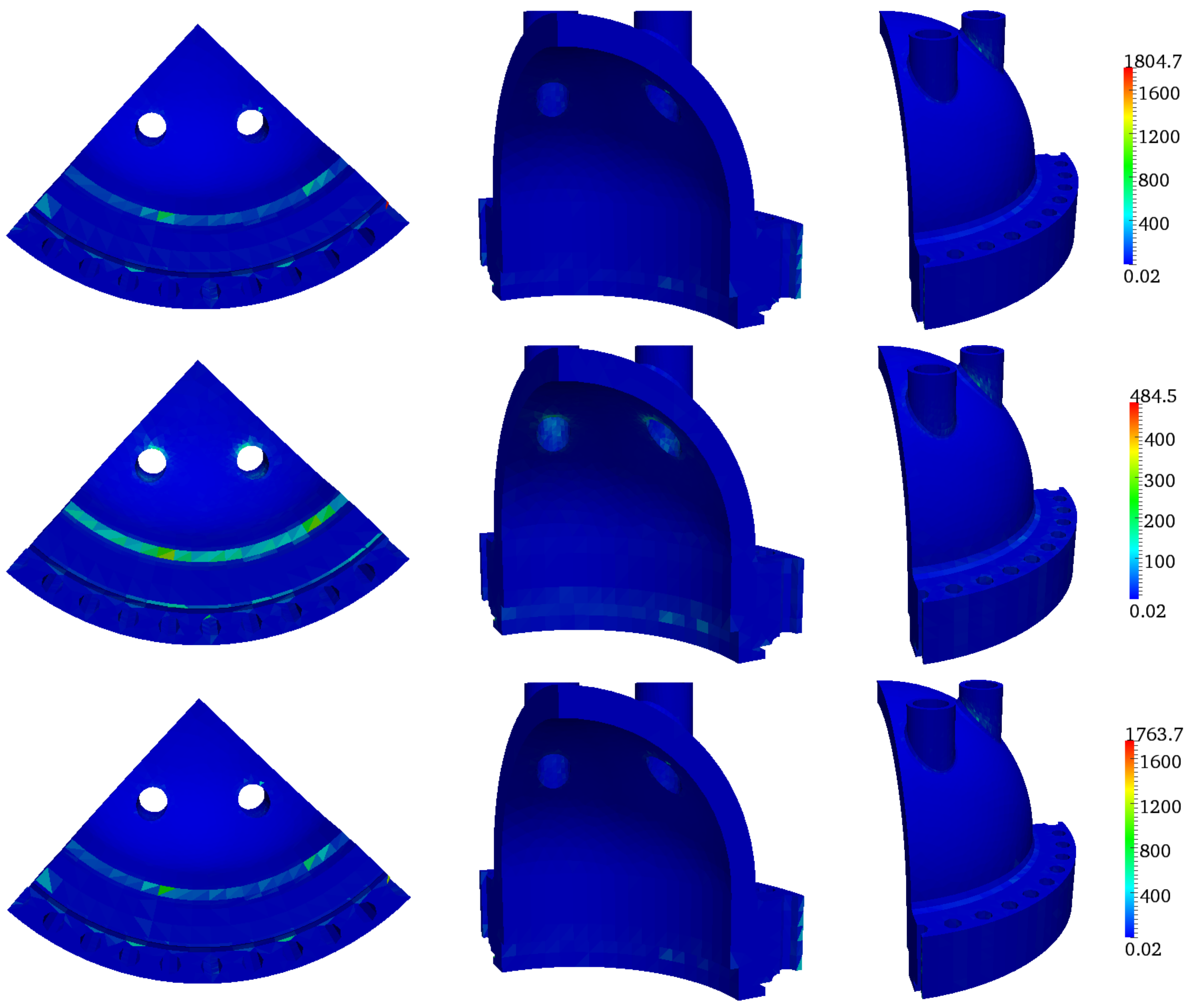}
\caption{Spatial distribution of local contributions to the density of the error estimates calculated using the EET (top), the SPET (middle), and the EESPT (bottom).}\label{fig1:distribution_estimators_tete_reacteur_3D}
\end{figure}

\begin{figure}
\centering\includegraphics[scale = 0.4]{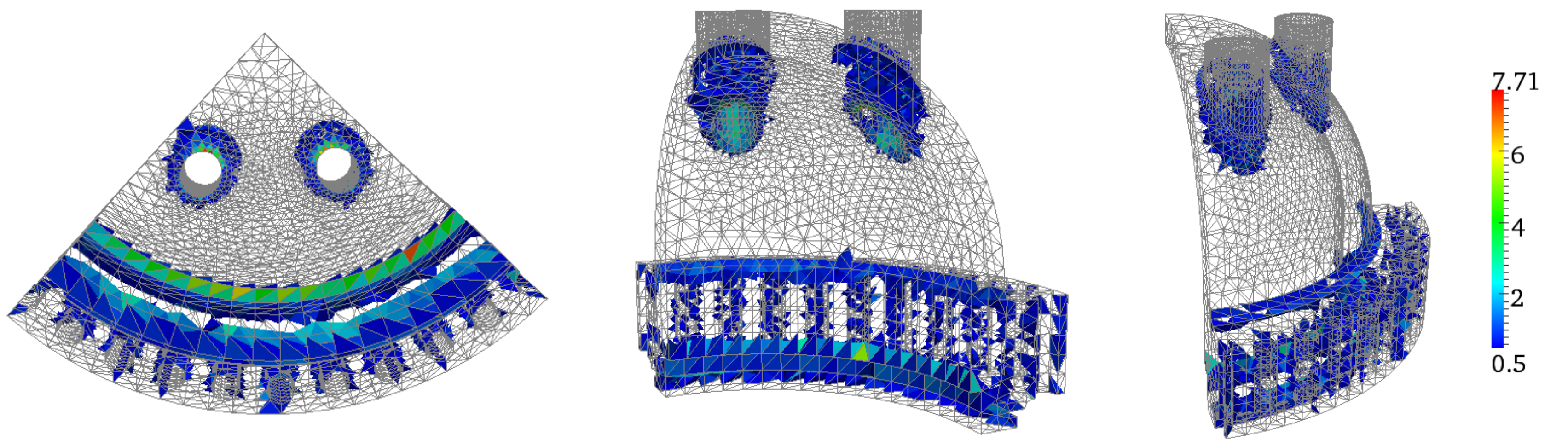}
\caption{Spatial distribution of relevant local contributions to the density of the energy norm of the reference error.}\label{fig1:distribution_discretization_error_tete_reacteur_3D_zoom}
\end{figure}

\begin{figure}
\centering\includegraphics[scale = 0.4]{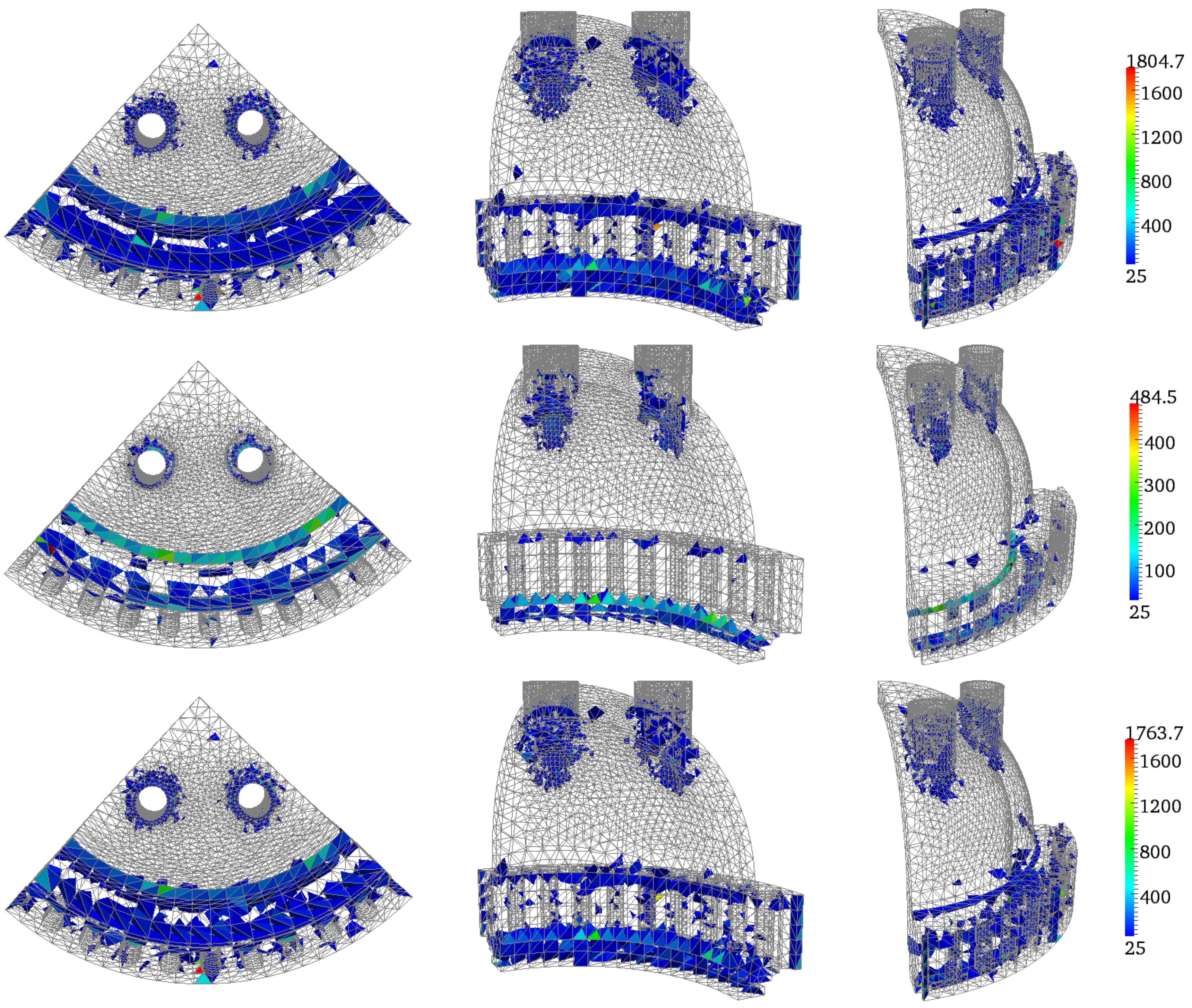}
\caption{Spatial distribution of relevant local contributions to the density of the error estimates calculated using the EET (top), the SPET (middle), and the EESPT (bottom).}\label{fig1:distribution_estimators_tete_reacteur_3D_zoom}
\end{figure}

The elementary contributions to the density of the energy norm of the reference error and that of the different error estimates are shown in \Figs{fig1:distribution_discretization_error_tete_reacteur_3D} and \ref{fig1:distribution_estimators_tete_reacteur_3D}, respectively. The higher contributions to the error estimates obtained using the EET and EESPT are concentrated in the ill-shaped elements. \Figs{fig1:distribution_discretization_error_tete_reacteur_3D_zoom} and \ref{fig1:distribution_estimators_tete_reacteur_3D_zoom} depict only the elements whose contribution to the energy norm of the reference error and to the error estimates is relevant. One can see that the error is localized not only in the bottom region of the standpipes but also near the area connecting the flange and the dome.

\newpage

\newpage

\newpage

\newpage
Once again, the local contributions to the effectivity index are computed with respect to the reference error. The local effectivity indices range between $1.163$ and $58.900$ for the EET; $1.074$ and $22.756$ for the SPET; $1.201$ and $51.274$ for the EESPT. One can see that all the elements have a local effectivity index superior to $1$ for the three techniques. Nevertheless, it is worthy noticing that the introduced global error estimators could not directly be used in a goal-oriented analysis, especially the EET and EESPT, because the maximum local contribution to the error estimates drastically overestimates the corresponding contribution to the energy norm of the reference error.

\subsubsection{Influence of the choice of the cost function involved in the EET and the EESPT}\label{5.3.2}

The study of the influence of the various cost functions on the performance of the two error estimators EET and EESPT has been investigated; results are summarized in \Tab{table1:influence_cost_function_closure_head}. They tend to show that the global effectivity index experiences a moderate decrease by $4 \%$ for the EET and a slight increase by $3 \%$ for the EESPT between the cost functions $J_0$ and $J_2$. The use of the cost function $J_2$ enables to obtain a slightly better effectivity index compared to the cost function $J_1$ for both estimators. As for the previous three-dimensional case, cost functions $J_0$ and $J_2$ provide a more relevant effectivity index for the EET and the EESPT, respectively.
\begin{table}
\caption{Influence of the choice of the cost function involved in the EET and the EESPT on the quality of the error estimators.}
\centering\small
\begin{tabular}{l p{0.2cm} c c c p{0.2cm} c c c}
\toprule
 & & \multicolumn{3}{c}{Estimate \ $\theta$} & & \multicolumn{3}{c}{Effectivity index \ $\eta$} \\
\ms \cline{3-5} \cline{7-9} \ms
Methods & & $J_0$ & $J_1$ &$J_2$ & & $J_0$ & $J_1$ & $J_2$ \\
\midrule
EET & & $812.771$ & $787.989$ & $782.711$ & & $6.001$ & $5.818$ & $5.779$ \\
EESPT & & $787.983$ & $815.111$ & $809.267$ & & $5.818$ & $6.019$ & $5.975$ \\
\bottomrule
\end{tabular}
\label{table1:influence_cost_function_closure_head}
\end{table}

\tabcap{18 cm}
\begin{table}
\caption{Influence of the choice of the cost function involved in the EET and the EESPT on the bounds of the local effectivity indices.}
\centering\small
\begin{tabular}{l p{0.2cm} c c c}
\toprule
 & & \multicolumn{3}{c}{Bounds of the local effectivity indices \ $\eta_E$} \\ 
\ms \cline{3-5} \ms
 Methods & & $J_0$ & $J_1$ &$J_2$ \\
\midrule
EET & & $1.163$ - $58.900$ & $1.201$ - $51.275$ & $1.095$ - $59.823$ \\
EESPT & & $1.201$ - $51.274$ & $0.948$ - $57.349$ & $0.878$ - $62.153$ \\
\bottomrule
\end{tabular}
\label{table1:influence_cost_function_i_eff_closure_head}
\end{table}

The bounds of the local effectivity indices obtained using each cost function are given in \Tab{table1:influence_cost_function_i_eff_closure_head} for both EET and EESPT methods. The analysis of the data reveals that the two global estimators are not well-suited to direct goal-oriented analysis for all the cost functions considered.

\newpage
\subsection{Hub of main rotor}\label{5.4}

\begin{figure}
\centering\includegraphics[scale = 0.4]{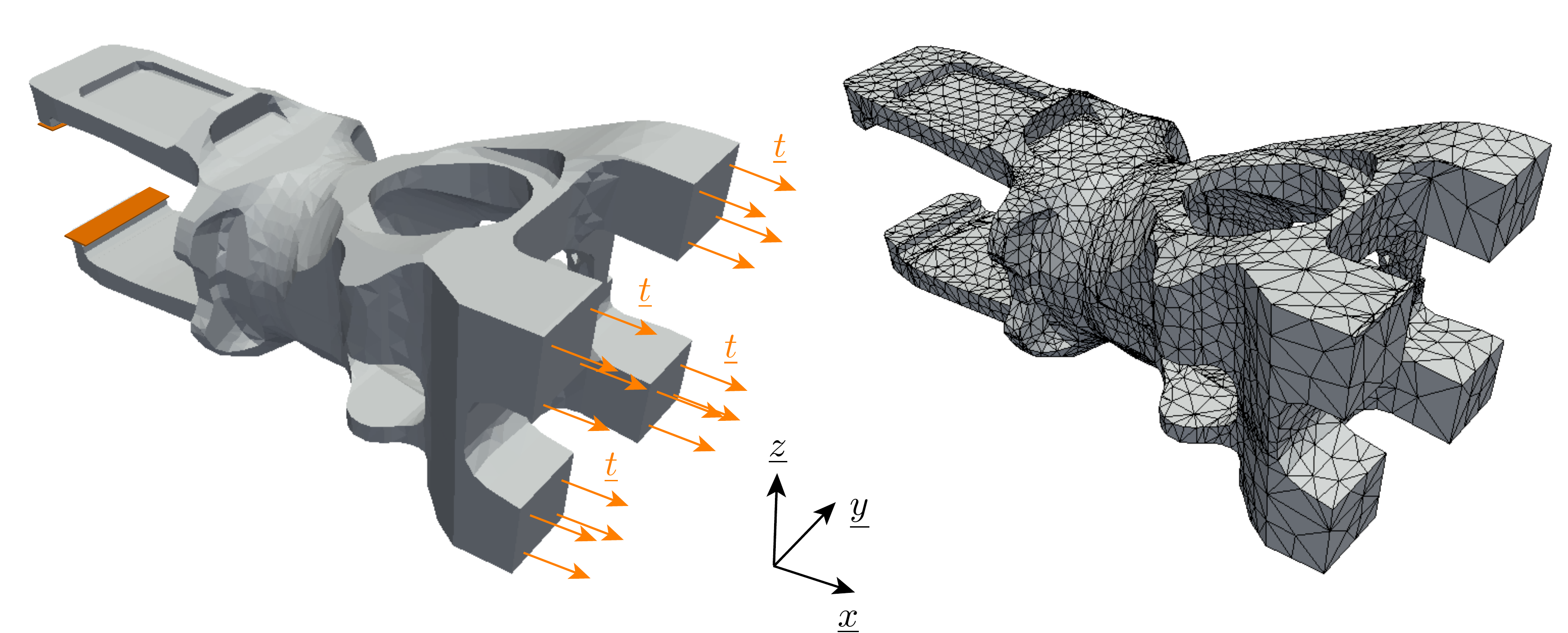}
\caption{Hub model problem (left) and associated finite element mesh (right). Orange plans represent clamped boundary conditions.}\label{fig1:manchon_helico_3D_geometry_mesh}
\end{figure}

A part of the NH90 helicopter from the Eurocopter company is considered. The structure is the hub of the main rotor which is used as a coupling sleeve between the helicopter frame and the rotor system. The structure is clamped at one end and subjected to a unit traction force density $\und{t}$, normal to the surface, on the other end. Let us notice that the loading plan is not exactly orthogonal to the main axis of the structure. The geometry and mesh considered, made of $19 \, 778$ linear tetrahedral elements and $5 \, 898$ nodes (\ie $17 \, 694$ d.o.f.), are shown in \Fig{fig1:manchon_helico_3D_geometry_mesh}. The reference mesh is built up by splitting each tetrahedron into $64$ tetrahedra. Therefore, it contains $1 \, 265 \, 792$ linear tetrahedral elements and $250 \, 274$ nodes (\ie $750 \, 822$ d.o.f.). One can notice that the FE mesh seems to be fairly distorted and, therefore, contains very ill-shaped elements.

\subsubsection{Comparison of the three error estimators}\label{5.4.1}

\begin{figure}
\centering\includegraphics[scale = 0.4]{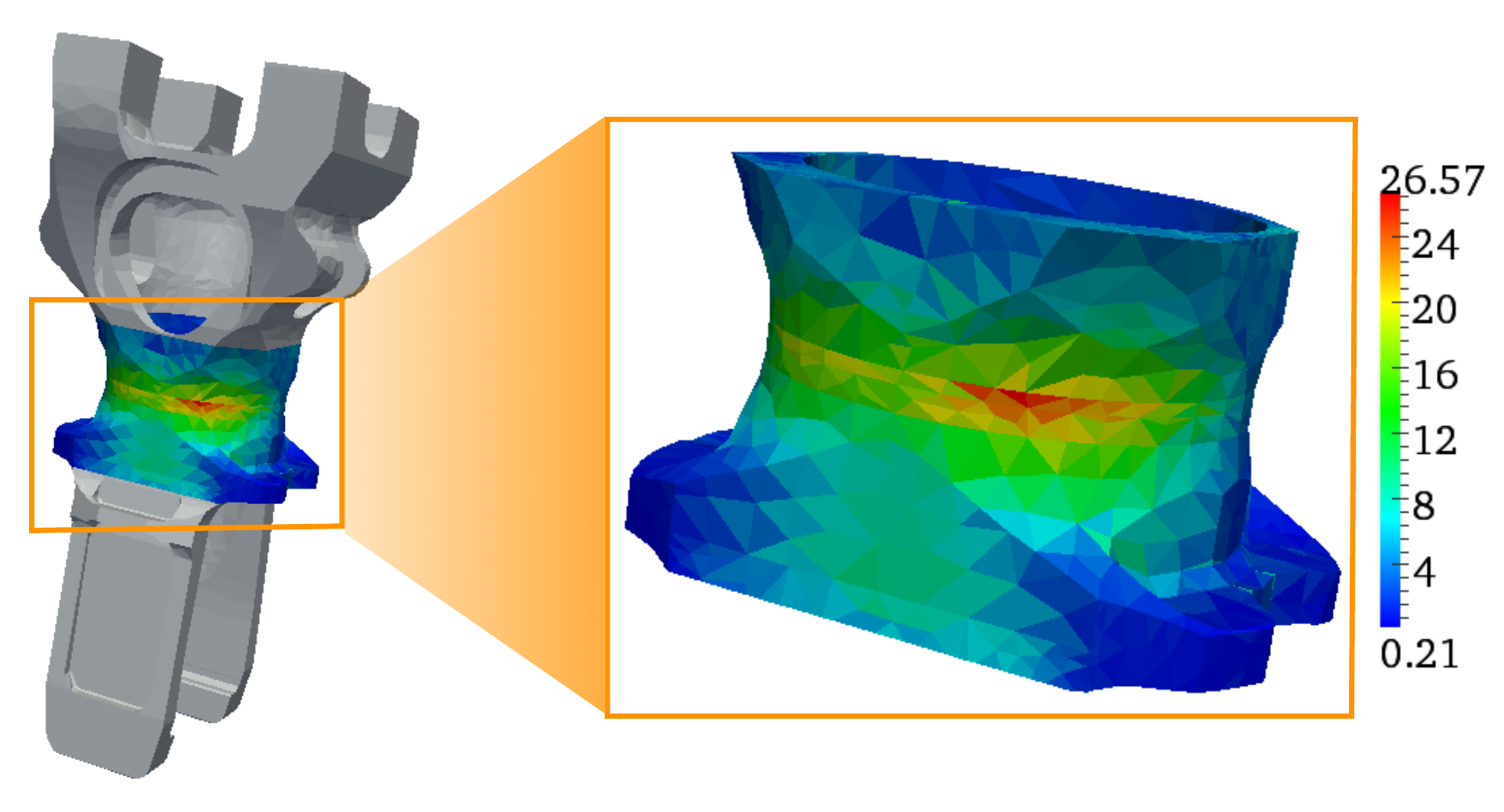}
\caption{Magnitude of the FE stress field.}\label{fig1:sigma_manchon_helico_3D}
\end{figure}

\begin{figure}
\centering\includegraphics[scale = 0.4]{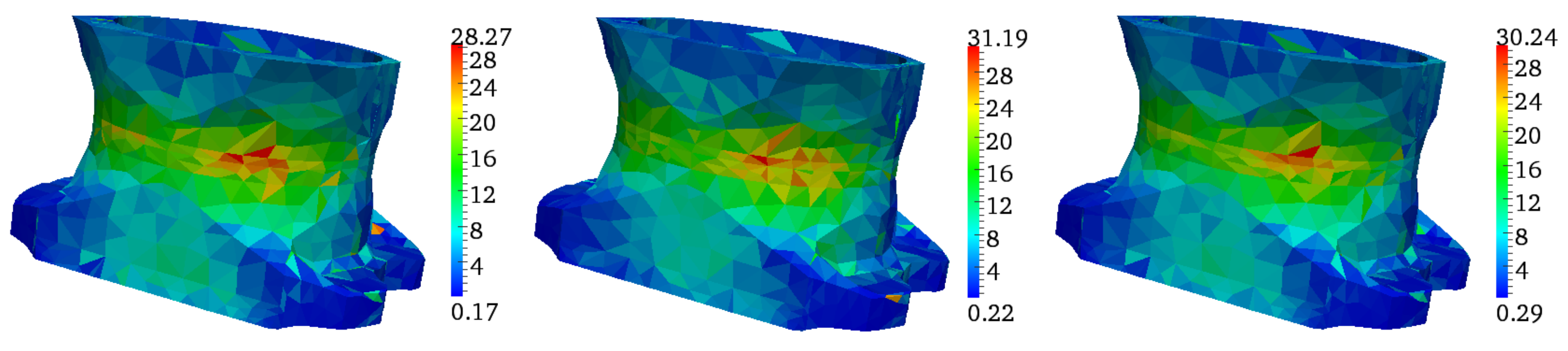}
\caption{Magnitude of the admissible stress field calculated using the EET (left), the SPET (center), and the EESPT (right).}\label{fig1:sigma_hat_manchon_helico_3D}
\end{figure}

The cost function $J_0$ has been used for the local minimization step. The highest stress region corresponds to the clamped surface, which is not a design zone. Conversely, the selected region in \Fig{fig1:sigma_manchon_helico_3D} plays an essential role in design purposes and engineering interest. The FE stress field in the selected region is depicted in \Fig{fig1:sigma_manchon_helico_3D} and the admissible stress fields obtained from the three techniques are displayed in \Fig{fig1:sigma_hat_manchon_helico_3D}.

The exact value of the energy norm of the reference error has been directly calculated from the reference solution:
\begin{align}
\lnorm{\und{e}_h}_{u, \Om} = \sqrt{\lnorm{\uu}^2_{u, \Om} - \lnorm{\uu_h}^2_{u, \Om}} \simeq 3852.53,
\end{align}
and required a computational cost of about $1$ s, whereas the computational cost necessary for calculating the local contributions to $\lnorm{\und{e}_h}_{u, \Om}$ reaches more than $14$ hours.

The estimated error assessment has been performed and the corresponding results are shown in \Tab{table1:comparison_hub_of_main_rotor}. The effectivity indices obtained using the EET and the EESPT are very high due to the high distortion of the FE mesh. Indeed, a large number of ill-shaped elements has a high contribution to the error. The error in the geometry may also be playing some role. Therefore, those error estimators are very sensitive to the bad quality of the FE mesh.

\begin{table}
\caption{Comparison of the error estimators given by the EET, the SPET, and the EESPT.}
\centering\small
\begin{tabular}{l c c c c}
\toprule
Methods & Estimate \ $\theta$ & Effectivity index \ $\eta$ & Normalized CPU time \\ 
\midrule
EET & $58 \, 061$ & $15.071$ & $1.000$ \\ 
SPET & $18 \, 948$ & $4.918$ & $4.828$ \\ 
EESPT & $42 \, 078$ & $10.922$ & $1.007$ \\ 
\bottomrule
\end{tabular}
\label{table1:comparison_hub_of_main_rotor}
\end{table}

\begin{figure}
\centering\includegraphics[scale = 0.4]{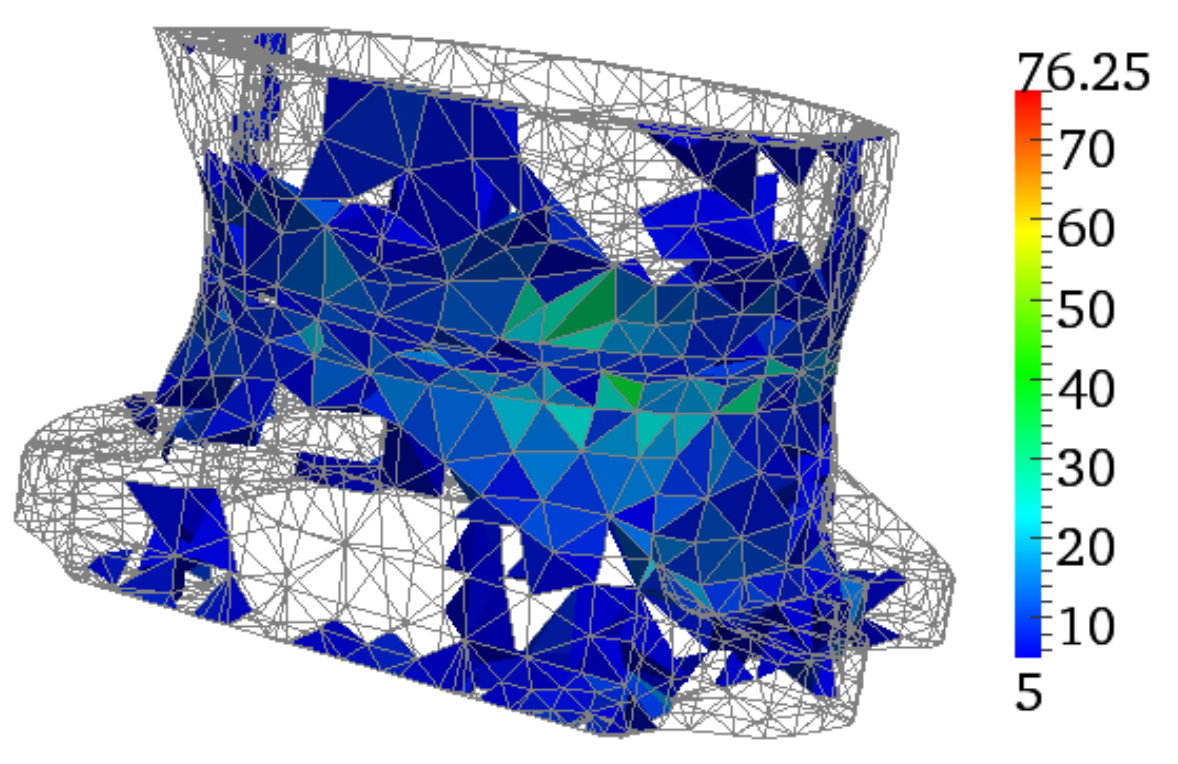}
\caption{Spatial distribution of relevant local contributions to the density of the energy norm of the reference error.}\label{fig1:distribution_discretization_error_manchon_helico_3D}
\end{figure}

\begin{figure}
\centering\includegraphics[scale = 0.38]{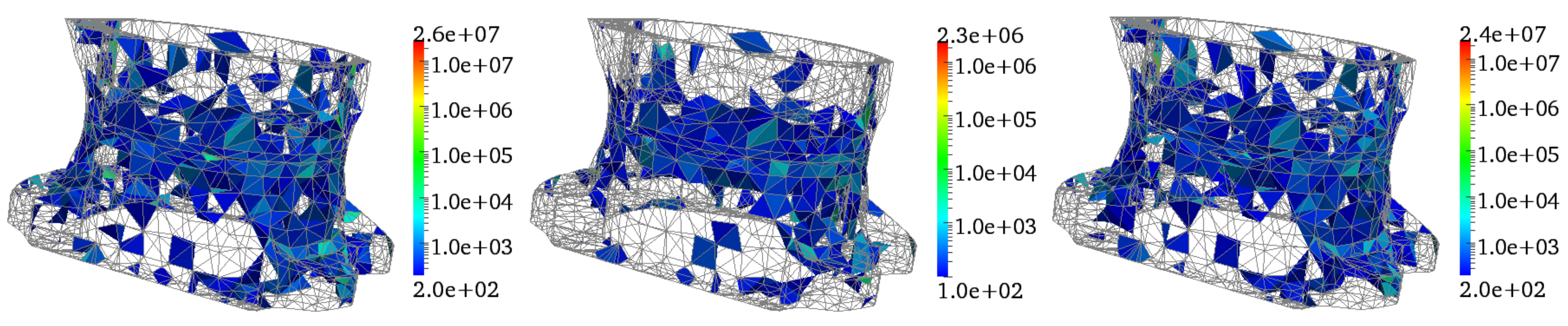}
\caption{Spatial distribution of relevant local contributions to the density of the error estimates calculated using the EET (left), the SPET (center), and the EESPT (right).}\label{fig1:distribution_estimators_manchon_helico_3D}
\end{figure}

The elementary contributions to the density of the energy norm of the reference error and that of the different error estimates are shown in \Figs{fig1:distribution_discretization_error_manchon_helico_3D} and \ref{fig1:distribution_estimators_manchon_helico_3D}, respectively. The contributions to the density of the various estimates are displayed on a log scale. Only elements located in zones where the contribution to the energy norm of the reference error and to the error estimates is significant are represented. Those maps reflect the same trend, namely a slightly higher contribution for the elements located in the highly-loaded region of the selected zone. The highest contributions are localized in the most distorted elements.


As for the previous cases, the local contributions to the effectivity index are computed with respect to the reference error. As the local effectivity indexes are very high for the most distorted elements, one distinguishes zones with very ill-shaped elements from the remainder of the FE mesh. The chosen quality measure used to evaluate tetrahedral quality is the radius ratio, which is the radius of the sphere circumscribed by a tetrahedron's four vertices divided by the radius of the inscribed sphere tangent to a tetrahedron's four edges. Most of the elements ($17 \, 773$ tetrahedra precisely) have a radius ratio less than $9$. For those elements, the local effectivity indices range between $0.25$ and $49.71$ for the EET; $0.51$ and $27.72$ for the SPET; $0.25$ and $47.15$ for the EESPT. On the other hand, for the most distorted elements ($2 \, 005$ tetrahedra exactly), which have a radius ratio varying between $9$ and $594.2$, the maximum value of the local effectivity indexes reaches $2529.25$ for the EET, $405.66$ for the SPET and $723.22$ for the EESPT. Overall, those numerical results indicate that the three estimators do not provide relevant local bounds in these regions.

\subsubsection{Influence of the choice of the cost function involved in the EET and the EESPT}\label{5.4.2}

Results are similar to that of the previous cases. \Tab{table1:influence_cost_function_hub_of_main_rotor} reports the estimates and the effectivity indexes obtained using the EET and the EESPT for the three cost functions we considered. The results show that the effectivity indexes are quite similar, except for EET and cost function $J_0$ whose corresponding effectivity index is fairly higher. 

\begin{table}
\caption{Influence of the choice of the cost function involved in the EET and the EESPT on the quality of the error estimators.}
\centering\small
\begin{tabular}{l p{0.2cm} c c c p{0.2cm} c c c}
\toprule
 & & \multicolumn{3}{c}{Estimate \ $\theta$} & & \multicolumn{3}{c}{Effectivity index \ $\eta$} \\
\ms \cline{3-5} \cline{7-9} \ms
Methods & & $J_0$ & $J_1$ &$J_2$ & & $J_0$ & $J_1$ & $J_2$ \\
\midrule
EET & & $58 \, 061$ & $42 \, 078$ & $42 \, 995$ & & $15.071$ & $10.922$ & $11.160$ \\
EESPT & & $42 \, 078$ & $42 \, 985$ & $43 \, 574$ & & $10.922$ & $11.158$ & $11.311$ \\
\bottomrule
\end{tabular}
\label{table1:influence_cost_function_hub_of_main_rotor}
\end{table}

\newpage
\section{CONCLUSION AND PROSPECTS}\label{6}

The three techniques, namely the EET, the SPET and the EESPT, have been implemented and analyzed in terms of quality of the error estimate, practical implementation and computational cost. Those methods have also been compared with respect to the energy norm of the reference error, i.e the global reference error.

Several two- and three-dimensional experiments have been carried out in order to compare the performance and computational cost of the three different approaches. The distribution of the local contributions to the global reference error has been accurately estimated for all the numerical experiments. Indeed, the three methods yield estimates which are guaranteed and sharp upper bounds of the energy norm of the reference error, even though the SPET appears to give a superior accuracy of estimation than that achievable by the EET and EESPT. Besides, the SPET seems to be convenient for implementation compared to the EET, and the EESPT to a lesser extent. However, the EET and the EESPT offer lower computational costs compared to the SPET, especially for three-dimensional cases. Thus, the EESPT may overcome the practical difficulties involved by the two other methods. Overall, the EESPT seems to be a good compromise in terms of quality of the error estimation, practical implementation and computational cost.

The development of an enhanced version of the EESPT method using a weak prolongation condition and resulting in local minimization of the complementary energy will be addressed in a forthcoming work; it is inspired from previous works of \cite{Flo02bis}. Future research will also focus on the use of EESPT method for goal-oriented error estimation associated to CRE.

\bibliographystyle{unsrt}

\begin{thebibliography}{10}

\bibitem{Lad75}
P.~Ladev{\`e}ze.
\newblock {\em {Comparaison de mod\`eles de milieux continus}}.
\newblock PhD thesis, Universit\'e Pierre et Marie Curie, 1975.

\bibitem{Lad83}
P.~Ladev{\`e}ze and D.~Leguillon.
\newblock {Error estimate procedure in the finite element method and
  application}.
\newblock {\em SIAM Journal of Numerical Analysis}, 20(3):485--509, 1983.

\bibitem{Bab78a}
I.~Babu\v{s}ka and W.~C. Rheinboldt.
\newblock {Error estimates for adaptative finite element computation}.
\newblock {\em SIAM, Journal of Numerical Analysis}, 15(4):736--754, 1978.

\bibitem{Zie87}
O.~C. Zienkiewicz and J.~Z. Zhu.
\newblock {A simple error estimator and adaptive procedure for practical
  engineerng analysis}.
\newblock {\em International Journal for Numerical Methods in Engineering},
  24(2):337--357, 1987.

\bibitem{Ver96}
R.~Verf{\"u}rth.
\newblock {\em {A Review of A Posteriori Error Estimation and Adaptive
  Mesh-refinement Techniques}}.
\newblock Wiley-Teubner, Stuttgart, 1996.

\bibitem{Bab01}
I.~Babu\v{s}ka and T.~Strouboulis.
\newblock {\em {The finite element method and its reliability}}.
\newblock Oxford University Press, Oxford, 2001.

\bibitem{Ste03}
E.~Stein, E.~Ramm, E.~Rank, R.~Rannacher, K.~Schweizerhof, E.~Stein,
  W.~Wendland, G.~Wittum, P.~Wriggers, and W.~Wunderlich.
\newblock {\em {Error-controlled Adaptative Finite Elements in Solid
  Mechanics}}.
\newblock J. Wiley, New York, 2003.

\bibitem{Lad04}
P.~Ladev{\`e}ze and J.~P. Pelle.
\newblock {\em {Mastering Calculations in Linear and Nonlinear Mechanics}}.
\newblock Springer, New York, 2004.

\bibitem{Par97}
M.~Paraschivoiu, J.~Peraire, and A.~T. Patera.
\newblock {A posteriori finite element bounds for linear-functional outputs of
  elliptic partial differential equations}.
\newblock {\em Computer Methods in Applied Mechanics and Engineering},
  150(1-4):289--312, 1997.

\bibitem{Ran97}
R.~Rannacher and F.-T. Suttmeier.
\newblock {A feed-back approach to error control in finite element methods:
  application to linear elasticity}.
\newblock {\em Computational Mechanics}, 19:434--446, 1997.
\newblock 10.1007/s004660050191.

\bibitem{Cir98}
F.~Cirak and E.~Ramm.
\newblock {A posteriori error estimation and adaptivity for linear elasticity
  using the reciprocal theorem}.
\newblock {\em Computer Methods in Applied Mechanics and Engineering},
  156:351--362, 1998.

\bibitem{Per98}
J.~Peraire and A.~T. Patera.
\newblock {Bounds for linear-functional outputs of coercive partial
  differential equations: Local indicators and adaptive refinement}.
\newblock In P.~Ladev{\`e}ze and J.T. Oden, editors, {\em Advances in Adaptive
  Computational Methods in Mechanics}, volume~47 of {\em Studies in Applied
  Mechanics}, pages 199--216. Elsevier, 1998.

\bibitem{Lad99a}
P.~Ladev{\`e}ze, P.~Rougeot, P.~Blanchard, and J.~P. Moreau.
\newblock {Local error estimators for finite element linear analysis}.
\newblock {\em Computer Methods in Applied Mechanics and Engineering},
  176(1-4):231--246, 1999.

\bibitem{Pru99}
S.~Prudhomme and J.~T. Oden.
\newblock {On goal-oriented error estimation for elliptic problems: application
  to the control of pointwise errors}.
\newblock {\em Computer Methods in Applied Mechanics and Engineering},
  176(1-4):313--331, 1999.

\bibitem{Str00}
T.~Strouboulis, I.~Babu\v{s}ka, D.~K. Datta, K.~Copps, and S.~K. Gangaraj.
\newblock {A posteriori estimation and adaptive control of the error in the
  quantity of interest. Part I: A posteriori estimation of the error in the von
  Mises stress and the stress intensity factor}.
\newblock {\em Computer Methods in Applied Mechanics and Engineering},
  181(1-3):261--294, 2000.

\bibitem{Bec01}
R.~Becker and R.~Rannacher.
\newblock {An optimal control approach to a posteriori error estimation in
  finite element methods}.
\newblock {\em Acta Numerica, A. Isereles (ed.), Cambridge University Press},
  10:1--120, 2001.

\bibitem{Wib06}
N.~E. Wiberg and P.~D{\'\i}ez.
\newblock {Special Issue}.
\newblock {\em Computer Methods in Applied Mechanics and Engineering},
  195:4--6, 2006.

\bibitem{Lad08}
Pierre Ladev{\`e}ze.
\newblock {Strict upper error bounds on computed outputs of interest in
  computational structural mechanics}.
\newblock {\em Computational Mechanics}, 42(2):271--286, 2008.

\bibitem{Cha07}
L.~Chamoin and Pierre Ladev{\`e}ze.
\newblock {Bounds on history-dependent or independent local quantities in
  viscoelasticity problems solved by approximate methods}.
\newblock {\em International Journal for Numerical Methods in Engineering},
  71(12):1387--1411, 2007.

\bibitem{Cha08}
L.~Chamoin and Pierre Ladev{\`e}ze.
\newblock {A non-intrusive method for the calculation of strict and efficient
  bounds of calculated outputs of interest in linear viscoelasticity problems}.
\newblock {\em Computer Methods in Applied Mechanics and Engineering},
  197(9-12):994--1014, 2008.

\bibitem{Pan10}
J.~Panetier, P.~Ladev{\`e}ze, and L.~Chamoin.
\newblock {Strict and effective bounds in goal-oriented error estimation
  applied to fracture mechanics problems solved with XFEM}.
\newblock {\em International Journal for Numerical Methods in Engineering},
  81(6):671--700, 2010.

\bibitem{Fra65}
B.~Fraeijs de~Veubeke.
\newblock {Displacement and equilibrium models in the finite element method by
  B. Fraeijs de Veubeke, Chapter 9, Pages 145--197 of Stress Analysis, Edited
  by O. C. Zienkiewicz and G. S. Holister, Published by John Wiley \& Sons,
  1965}.
\newblock {\em International Journal for Numerical Methods in Engineering},
  52(3):287--342, 1965.

\bibitem{Kem03}
M.~Kempeneers, P.~Beckers, J.~P.~Moitinho de~Almeida, and O.~J. B.~A. Pereira.
\newblock {Mod\`eles \'equilibre pour l'analyse duale}.
\newblock {\em Revue Europ\'eenne des \'Elements Finis}, 12(6):737--760, 2003.

\bibitem{Ain93bis}
M.~Ainsworth and J.~T. Oden.
\newblock {A posteriori error estimators for second order elliptic systems:
  Part 2. An optimal order process for calculating self-equilibrating fluxes}.
\newblock {\em Computers \& Mathematics with Applications}, 26(9):75--87, 1993.

\bibitem{Lad96}
P.~Ladev{\`e}ze and E.~A.~W. Maunder.
\newblock {A general method for recovering equilibrating element tractions}.
\newblock {\em Computer Methods in Applied Mechanics and Engineering},
  137(2):111--151, 1996.

\bibitem{Lad97}
P.~Ladev{\`e}ze and P.~Rougeot.
\newblock {New advances on a posteriori error on constitutive relation in
  finite element analysis}.
\newblock {\em Computer Methods in Applied Mechanics and Engineering},
  150(1-4):239--249, 1997.

\bibitem{Flo02}
E.~Florentin, L.~Gallimard, and J.~P. Pelle.
\newblock {Evaluation of the local quality of stresses in 3D finite element
  analysis}.
\newblock {\em Computer Methods in Applied Mechanics and Engineering},
  191(39-40):4441--4457, 2002.

\bibitem{Mac00}
Luc Machiels, Yvon Maday, and Anthony~T. Patera.
\newblock {A ``flux-free'' nodal Neumann subproblem approach to output bounds
  for partial differential equations}.
\newblock {\em Comptes Rendus de l'Acad\'emie des Sciences - Series I -
  Mathematics}, 330(3):249--254, 2000.

\bibitem{Car00}
Carsten Carstensen and Stefan~A. Funken.
\newblock {Fully Reliable Localized Error Control in the FEM}.
\newblock {\em SIAM J. Sci. Comput.}, 21(4):1465--1484, 2000.

\bibitem{Mor03}
P.~Morin, R.~H. Nochetto, and K.~G. Siebert.
\newblock {Local problems on stars: a posteriori error estimators, convergence,
  and performance}.
\newblock {\em Mathematics of Computation}, 72(243):1067--1097, 2003.

\bibitem{Pru04}
S.~Prudhomme, F.~Nobile, L.~Chamoin, and J.~T. Oden.
\newblock {Analysis of a subdomain-based error estimator for finite element
  approximations of elliptic problems}.
\newblock {\em Numerical Methods for Partial Differential Equations},
  20(2):165--192, 2004.

\bibitem{Par06}
N.~Par{\'e}s, P.~D{\'\i}ez, and A.~Huerta.
\newblock {Subdomain-based flux-free a posteriori error estimators}.
\newblock {\em Computer Methods in Applied Mechanics and Engineering},
  195(4-6):297--323, 2006.

\bibitem{Par09}
N.~Par{\'e}s, H.~Santos, and P.~D{\'\i}ez.
\newblock {Guaranteed energy error bounds for the Poisson equation using a
  flux-free approach: Solving the local problems in subdomains}.
\newblock {\em International Journal for Numerical Methods in Engineering},
  79(10):1203--1244, 2009.

\bibitem{Moi09}
J.~P. Moitinho~de Almeida and E.~A.~W. Maunder.
\newblock {Recovery of equilibrium on star patches using a partition of unity
  technique}.
\newblock {\em International Journal for Numerical Methods in Engineering},
  79(12):1493--1516, 2009.

\bibitem{Cot09}
R.~Cottereau, P.~D\'iez, and A.~Huerta.
\newblock {Strict error bounds for linear solid mechanics problems using a
  subdomain-based flux-free method}.
\newblock {\em Computational Mechanics}, 44(4):533--547, 2009.

\bibitem{Gal09}
L.~Gallimard.
\newblock {A constitutive relation error estimator based on traction-free
  recovery of the equilibrated stress}.
\newblock {\em International Journal for Numerical Methods in Engineering},
  78(4):460--482, 2009.

\bibitem{Lad08bis}
P.~Ladev{\`e}ze.
\newblock {V\'erification des calculs \'el\'ements finis : une nouvelle
  technique non intrusive pour le calcul des contraintes admissibles}.
\newblock Technical Report 270, LMT-Cachan, 2008.

\bibitem{Lad10bis}
P.~Ladev{\`e}ze, L.~Chamoin, and E.~Florentin.
\newblock {A new non-intrusive technique for the construction of admissible
  stress fields in model verification}.
\newblock {\em Computer Methods in Applied Mechanics and Engineering},
  199:766--777, 2010.

\bibitem{Lad92}
P.~Ladev{\`e}ze, P.~Marin, J.~P. Pelle, and J.~L. Gastine.
\newblock {Accuracy and optimal meshes in finite element computation for nearly
  incompressible materials}.
\newblock {\em Computer Methods in Applied Mechanics and Engineering},
  94(3):303--315, 1992.

\bibitem{Pra47}
W.~Prager and J.~L. Synge.
\newblock {Approximations in elasticity based on the concept of functions
  spaces}.
\newblock {\em Quarterly of Applied Mathematics}, 5:261--269, 1947.

\bibitem{Moi06}
J.~P.~Moitinho de~Almeida and O.~J. B.~Almeida Pereira.
\newblock {Upper bounds of the error in local quantities using equilibrated and
  compatible finite element solutions for linear elastic problems}.
\newblock {\em Computer Methods in Applied Mechanics and Engineering},
  195(4-6):279--296, 2006.
\newblock Adaptive Modeling and Simulation.

\bibitem{Bab94}
I.~Babu\v{s}ka, T.~Strouboulis, C.~S. Upadhyay, S.~K. Gangaraj, and K.~Copps.
\newblock {Validation of a posteriori error estimators by numerical approach}.
\newblock {\em International journal for numerical methods in engineering},
  37(7):1073--1123, 1994.

\bibitem{Coo99}
P.~Coorevits, J.~P. Dumeau, and J.~P. Pelle.
\newblock {Control of analyses with isoparametric elements in both 2D and 3D}.
\newblock {\em International journal for numerical methods in engineering},
  46(2):157--176, 1999.

\bibitem{Samtech}
{www.samtech.com}.

\bibitem{Flo02bis}
E.~Florentin.
\newblock {\em {Sur l'\'evaluation de la qualit\'e locale des contraintes
  \'el\'ements finis en \'elasticit{\'e} tridimensionnelle.}}
\newblock PhD thesis, \'Ecole Normale Sup\'erieure de Cachan, 2002.

\end{thebibliography}

\end{document}